\documentclass[twoside,openright,11pt,a4paper,epsf]{article}
\usepackage{latexsym}
\usepackage[T1]{fontenc}
\usepackage{amsmath}
\usepackage{amsthm}
\usepackage[french]{babel}
\usepackage{amsfonts}
\usepackage{amssymb}
\usepackage{vmargin}
\usepackage{fancyheadings}
\usepackage{epsfig}
\usepackage[all]{xy}

\newcommand{\color}[6]{}
\newcommand{\Q}{\mathbb{Q}}
\newcommand{\R}{\mathbb{R}}
\newcommand{\N}{\mathbb{N}}

\newcommand{\C}{\mathbb{C}}
\newcommand{\D}{\mathbb{D}}
\newcommand{\B}{\mathbb{B}}
\renewcommand{\P}{\mathbb{P}}
\renewcommand{\S}{\mathbb{S}}

\newcommand{\grad}{\vec{\nabla}}

\newcommand{\priv}{\backslash}
\newcommand{\lra}{\longrightarrow}

\newcommand{\om}{\Omega}
\newcommand{\eps}{\epsilon}
\newcommand{\F}{\mathcal{F}}
\renewcommand{\L}{\mathcal{L}}
\newcommand{\psh}{$p.s.h$}
\newcommand{\wdt}[1]{\widetilde{#1}}
\newcommand{\cqfd}{\hfill $\square$ \vspace{0.1cm}\\ }
\newcommand{\sbull}{{\tiny $\bullet$ }}
\newcommand{\ds}{\displaystyle}
\newcommand{\spc}{\mathcal{SPC}}
\newcommand{\fpc}{\mathcal{FPC}}
\newcommand{\im}{\text{Im}\,}
\newcommand{\re}{\text{Re}\,}

\newtheorem{definition}{D\'efinition }[section]

\newtheorem{lemme}[definition]{Lemme}
\newtheorem*{question}{Question}
\newtheorem{theoreme}[definition]{Th\'eor\`eme}

\newtheorem{prop}[definition]{Proposition}
\newtheorem{rque}[definition]{Remarque}

\newcounter{faitmoins}

\newtheorem*{thmncgal}{Théorème 2}
\newtheorem*{thmrec}{Théorème 1}
\newtheorem*{thmdisque}{Théorème 3}

\title{Dynamique des applications holomorphes propres de domaines réguliers et
  problème de l'injectivité.}
\author{Emmanuel Opshtein.}
\begin{document}
\maketitle
\section*{Introduction}
En 1977, Alexander montra que les auto-applications
holomorphes propres de la boule euclidienne de $\C^k$ ($k\geq 2$) sont des
automorphismes \cite{alexander}. Peu avant, Pinchuk (\cite{pinchuk}, 1974)
avait établi un résultat analogue pour les domaines strictement pseudoconvexes
à bords lisses 
lorsque les applications se prolongent différentiablement au bord. La
découverte d'un phénomène tranchant aussi singulièrement avec  le cas
unidimensionnel souleva une question qui est demeurée ouverte : 
\begin{question}
Une auto-application holomorphe propre d'un domaine pseudoconvexe borné à bord
lisse dans $\C^k$ est-elle nécessairement un automorphisme lorsque $k>1$?
\end{question}
L'hypothèse de régularité du bord est évidemment nécessaire comme le montrent
l'exemple de $(w,z^2)$ sur le bidisque ou d'autres bassins d'attraction. Une
première série de travaux s'inscrivit dans la perspective ouverte par
Alexander. Ainsi, Pinchuk (\cite{pinchuk3}, 1978) montra comment ramener le cas
strictement pseudoconvexe à celui de la boule par un argument de dilatations de
coordonnées. Henkin et Tumanov (\cite{henktum}, 1982) puis Henkin et Novikov 
(\cite{henknov}, 1984) généralisèrent le résultat d'Alexander aux domaines
symétriques bornés.   
Les progrès accomplis dans les années 80 sur le prolongement des applications
holomorphes propres ont considérablement simplifié l'abord de cette question.
En particulier, le théorème d'Alexander devient limpide car les applications
se prolongent au voisinage de la frontière (\cite{bellreinhardt}, 1983).
Dans ce contexte, une nouvelle vague de travaux aborda le cas des domaines
faiblement pseudoconvexes.
Citons les plus marquants.
Bedford et Bell (\cite{bedbell}, 1982) ont résolu le cas où $b\om$ est analytique
réel. Diederich et Fornaess (\cite{diederichfornaess1}, 1982) celui où le lieu de
faible pseudoconvexité est petit (de mesure $2k-3$ dimensionnelle
nulle). D'autres réponses concernent des domaines présentant des symétries. 
Citons Berteloot et Pinchuk (\cite{bertpinchuk}, 1994) pour les domaines de
Reinhardt de $\C^2$ bornés complets et distincts du bidisque; Berteloot
(\cite{berteloot1}, 1998) pour les domaines de Reinhardt de $\C^k$ bornés à
frontière $\mathcal C^2$. 
Coupet, Pan et Sukhov (\cite{cps}, 2001) ont traité le cas des 
domaines cerclés de type fini de $\C^2$.
Enfin, dans un registre sensiblement différent, Tsaï (\cite{tsai}, 1993) retrouve et
prolonge les résultats de Henkin-Novikov. Le lecteur pourra consulter
l'article de revue de Forstneri\v{c} \cite{forstneric} pour d'autres
références. 

La plupart de ces travaux reposent sur une stratégie commune exploitant la
compétition entre deux phénomènes : l'ensemble critique croît avec l'itération
 mais sa trace au bord reste limitée par le lieu de faible pseudoconvexité. 
\newpage

Dans cet article, nous adoptons un nouveau point de vue axé sur une approche
plus dynamique. Nous montrons que l'entropie des prolongements au bord est
toujours portée par le lieu de faible pseudoconvexité et confrontons cette
contrainte à la richesse dynamique en présence de degré. 

Dans cette perspective, nous étudions le lien entre la dynamique d'une
application au bord et à l'intérieur d'un domaine. Cette dernière est
gouvernée par une dichotomie remarquable. Elle est soit non-récurrente,
ce qui signifie que toutes les limites sont à valeurs dans le bord, soit
récurrente auquel cas l'orbite de tout compact reste relativement
compacte. Dans le cas d'une dynamique récurrente, toutes les limites sont des
automorphismes d'une même variété holomorphe obtenue par rétraction du
domaine. Nous dirons que cette variété est le rétract de
l'application. \vspace{0.2cm}

Notre premier résultat montre que la régularité du bord est une obstruction à
la récurrence :
\begin{thmrec}
Soit $\om$ un domaine borné de $\C^k$ possédant une fonction \psh\ lisse et
définissante globale. Soit 
$f:\om\lra \om$ une application holomorphe propre dont la dynamique est
récurrente. Le rétract de $f$ ne peut être ni un point ni une surface de
Riemann. Si $k=2$ alors $f$ est un automorphisme.
\end{thmrec}
Soulignons que l'essentiel de la démonstration de ce théorème est valable en
toutes dimensions.

La démonstration de ce théorème s'adapte aux dynamiques
non-récurrentes en se restreignant à une classe de domaines possédant de
bonnes fonctions pics \psh : les domaines $LB$-réguliers. Cette classe englobe
les domaines de type 
fini de $\C^2$ et les domaines strictement géométriquement convexes.  Nous
montrons que l'entropie de $f_{|b\om}$ est portée par le lieu de faible
pseudoconvexité en contrôlant son ensemble non-errant $NW(f_{|b\om})$ :
\begin{thmncgal}
Soit $\om$ un domaine borné de $\C^k$ à bord lisse et $LB$-régulier dont le
bord contient au moins deux points de faible pseudoconvexité. Soit $f:\om\lra
\om$ holomorphe propre se prolongeant différentiablement au bord et dont la
dynamique est non-récurrente. Alors l'alternative suivante a lieu :
\begin{enumerate}
\item les limites de $f^n$ sur $\om$ et sur $\spc(b\om)$ sont des points
  identiques de $\fpc(b\om)$ et $NW(f_{|b\om})\subset \fpc(b\om)$ 
\item[ou]
\item $f^n$ converge localement uniformément sur $\om$ et sur $\spc(b\om)$ vers
  $a\in \spc(b\om)$ et $NW(f_{|b\om})\subset \fpc(b\om)\cup\{a\}$. 
\end{enumerate}

\end{thmncgal}
Enfin, nous faisons aboutir notre stratégie pour la classe des domaines
disqués de $\C^2$. Nous obtenons une généralisation  du résultat de
Coupet-Pan-Sukhov \cite{cps} :
\begin{thmdisque}
Soit $\om\Subset \C^2$ un domaine pseudoconvexe disqué, 
$LB$-régulier et dont la fonction de jauge est lisse. Alors toute application
holomorphe propre de $\om$ dans lui-même est un automorphisme de $\om$.  
\end{thmdisque}
La démonstration de ce théorème épouse parfaitement le point de vue développé
par cet article. Lorsque $f$ n'est pas injective on voit, en utilisant les
théorèmes 1 et 2 ainsi que la géométrie du domaine, que toute la dynamique de
$f_{|b\om}$ se concentre sur un unique cercle de faible pseudoconvexité. Des
estimations standards d'entropie topologique montrent que ceci est impossible. 
\paragraph{} L'article est organisé de la façon suivante. La
première partie concerne la dynamique des applications considérées. Nous y
précisons la dichotomie récurrence/non-récurrence et présentons certains exemples 
typiques. La deuxième partie contient des  préliminaires  techniques
standards. Dans la troisième partie, nous expliquons pourquoi le bord du
bassin d'un point attractif ne peut être lisse. Ce cas simple offre une
esquisse de démonstration pour les théorèmes 1 et 2.
Dans la quatrième partie, nous établissons des estimées de dérivées pour les
applications $CR$. Nous utilisons ces estimées dans la cinquième partie pour
mettre en évidence l'expansivité des  systèmes dynamiques récurrents. Nous en
déduisons immédiatement le théorème 1. Dans la sixième 
partie, nous adaptons les raisonnements précédents pour obtenir le théorème
2. Enfin, la septième partie traite  des domaines disqués et  contient la
preuve du théorème 3. 

 Je tiens à remercier Francois Berteloot pour les
suggestions mathématiques et stylistiques dont il m'a fait bénéficier.   

\paragraph{Notations :}
\begin{itemize}
\item $\spc(b\om)$ est l'ensemble des points de stricte pseudoconvexité de
  $b\om$. 
\item $\fpc(b\om)$ est l'ensemble des points de faible pseudoconvexité de
  $b\om$.  
\item Si $x\in \fpc(b\om)$, la composante connexe de   $x$ dans $\fpc(b\om)$
  est un compact de $b\om$ noté $\fpc(x)$   
\item Si $p\in b\om$, $\vec N(p)$ est la normale à $b\om$ en $p$, dirigée vers
  l'intérieur de $\om$.
\item Si $z_1,z_2\in \C^k$, $\D(z_1,z_2)$ est le disque de $z_1+\C\cdot
  (z_2-z_1)$ de diamètre $[z_1,z_2]$.
\item $\mathcal C(z_1,z_2)=b\D(z_1,z_2)$.
\item $\mathcal C^\alpha(z_1,z_2)$ est l'arc de cercle de
  $\mathcal C(z_1,z_2)=b\D(z_1,z_2)$ centré en $z_1$ d'ouverture angulaire
  $2\alpha$. 
\item $\D_r:=\{z\in \C,\; |z|<r\}$.
\item $d_{K_\om}$ (ou $d_K$ en l'absence d'ambiguïté) est la distance de
  Kobayashi dans $\om$.
\end{itemize}
\section{Dynamiques holomorphes dans un domaine.}
Comme le montre le résultat suivant, la dynamique d'une application holomorphe
d'un domaine dans lui-même est soumise à une alternative simple :
\begin{theoreme}\label{abate}
Si $\om$ est un domaine taut de $\C^k$ et $f$ est une application holomorphe
de $\om$ dans lui-même alors :
\begin{enumerate}
\item $d\big(f^n(z),b\om\big)\lra 0$ localement uniformément
\item[ou]
\item il existe une sous-variété holomorphe $M$ de $\om$ et, parmi les limites
  de $(f^n)_n$, une rétraction holomorphe $\rho:\om\lra M$ ($\rho\circ
  \rho=\rho$). De plus, $f_{|M}\in$Aut$\,M$ et toutes les limites de $(f^n)_n$
  sont de la forme   $\gamma \circ \rho$ où $\gamma\in$Aut$\,M$. 
\end{enumerate}
\end{theoreme}
Dans le premier cas nous dirons que la dynamique de $f$ est {\it
  non-récurrente}, dans le second qu'elle est {\it récurrente}.

Nous renvoyons au livre de M. Abate \cite{abate}  pour une démonstration de ce
résultat (théorèmes 2.1.29 et 2.4.3). Soulignons que les domaines
pseudoconvexes bornés à bord lisse sont taut.
\begin{rque}\label{pasaut}
Lorsque $\dim M<k$, on peut penser aux fibres de la rétraction comme
aux variétés stables de $f$; $f$ échange ces fibres ($\rho\circ f=f\circ
\rho$) et sa dynamique y est contractante. Comme les automorphismes sont des
isométries pour la distance de Kobayashi, ceci ne peut se produire que lorsque
$f$ est non-injective. Par contre, lorsque la dimension de $M$ est maximale, $M$
coïncide avec $\om$, $\rho$ avec l'identité et $f$ est un automorphisme de
$\om$.  
\end{rque}
Décrivons quelques exemples de dynamique récurrente.\vspace{0.2cm}\\ 
\textbf{1. Exemples de Lattès.}
Dans la situation extrême où le rétract est de dimension nulle, la suite  
$(f^n)_n$ converge localement uniformément vers un point de $\om$. Tel est le
cas des applications polynomiales homogènes et non-dégénérées de $\C^k$ car
elles possèdent un bassin d'attraction superattractif à l'origine. Parmi
celles-ci, un exemple particulièrement frappant, étudié par C. Dupont dans
\cite{tof}, est obtenu en relevant à $\C^k$ un endomorphisme de Lattès de
$\P^{k-1}$. Pour cet exemple, le bord du bassin superattractif est lisse et
sphérique en dehors de la trace d'un ensemble algébrique. \vspace{0.2cm}\\  
\textbf{2. Exemple de Yoccoz; persistence des variétés stables.}
L'application $(w,z)\longmapsto (w,wz+z^2)$ produit un rétract
($\D\times \{0\}$) de dimension $1$ dans $\C^2$ et induit une auto-application
propre sur le bassin d'attraction de celui-ci. Yoccoz montre  par un argument
élémentaire que, pour la plupart des $\xi\in b\D$, les variétés stables des
$(w,0)$ ne disparaissent pas quand $w$ tend 
radialement vers $\xi$ et produisent une trace analytique dans le bord. Il
établit ainsi très simplement l'existence de disques de Siegel
(\cite{yoccoz}).     

Ce phénomène de persistence des variétés stables est remarquable et on peut
se demander s'il est général. Ceci entraînerait que le bord d'un domaine
soutenant une dynamique récurrente de rétract non-trivial contient des
variétés analytiques. La famille que nous décrivons maintenant montre que ce
phénomène de persistence n'a pas toujours lieu. Elle est obtenue en perturbant
l'exemple de Yoccoz.\vspace{0.2cm}\\
\textbf{3. La famille $f_\alpha=(we^{i2\pi\alpha},wz+z^2)$; évanescence des
  variétés stables.}
Pour $\alpha\in \R$, définissons
$$ 
\begin{array}{l}
\begin{array}{rcccl}
 f_\alpha & : & \C^2 & \longrightarrow &\C^2\\
  & & (w,z) & \longmapsto & (we^{i2\pi
 \alpha},wz+z^2)
\end{array}\vspace{0.2
cm}\\
\text{et } B_\alpha=\{(w,z)\in \C^2\; |\; f_{\alpha,2}^n(w,z)\lra 0\} \text{ où
  } f^n_\alpha:=(f^n_{\alpha,1},f^n_{\alpha,2}).
\end{array}
$$
 Le domaine
$B_\alpha$ est pseudoconvexe et contenu dans $\D\times \D_2$. L'application
$f_\alpha$ induit une application propre de $B_\alpha$ dans lui-même. On
pourra observer que $B_\alpha$ contient le triangle $T:=\{|w|+|z|<1\}$ et que
$f_\alpha(T)\subset T$, si bien que $B_\alpha=\cup f_\alpha^{-n}(T)$. La
dynamique de $f_\alpha$ est récurrente sur $B_\alpha$, la rétraction est
donnée par $\rho(w,z)=(w,0)$ et $M=\D\times\{0\}$. Pour tout $w\in \D$,
nous noterons $\F_w$ la fibre de $\rho$ au dessus de $(w,0)$.

Pour $\alpha=0$, on retrouve l'exemple de Yoccoz et la fibre $\F_w$ coïncide
avec le bassin d'attraction de $z\mapsto wz+z^2$ à l'origine. 

Nous allons voir, par un argument de catégories de Baire, que ces fibres
disparaissent pour certaines valeurs de $\alpha$. Plus précisément, {\it le rayon
conforme} $r(\alpha,w)$ de $\F_w$ en $(w,0)$ tend vers $0$ lorsque $|w|$ tend
vers $1$. A cet effet, introduisons les fonctions intermédiaires suivantes.
$$
r(\alpha,w) :=\sup\left\{r\;|\; \exists \psi\in {\cal O}(\D_r,\{w\}\times \C)
 \textit{ t.q. } \left|
\begin{array}{l} 
i) \;\psi(0)=(w,0) \text{ et } |\psi'(0)|=1\\
ii)\; \psi(\D_r)\subset \overline{B_\alpha}
\end{array}\right.
\right\}.
$$
On définit de même $\wdt r(\alpha,w)$ en remplaçant ii) par une condition moins
contraignante :
$$
\wdt r(\alpha,w) :=\sup\left\{r\;|\; \exists \psi\in {\cal O}(\D_r,\{w\}\times \C)
 \textit{ t.q. }\left| 
\begin{array}{l} 
i) \;\psi(0)=(w,0) \text{ et } |\psi'(0)|=1\\
ii)\; f_\alpha^n(\psi(\D_r))\subset \overline{\D\times \D_2}\hspace{0.2cm} \forall n\in \N 
\end{array}\right.
\right\}.
$$
On pose ensuite
$$
\begin{array}{lll}
R(\alpha) & := & \ds\inf_{w\in b\D} r(\alpha,w)\\
\wdt R(\alpha) & :=& \ds\inf_{w\in b\D} \wdt r(\alpha,w).
\end{array}
$$
On vérifie facilement le lemme suivant : 
\begin{lemme}\label{scs}
  Les fonctions $r(\alpha,\cdot)$, $\wdt r$ et $\wdt R$ sont semi-continues
  supérieurement sur $\overline{\D}$, $\R\times \overline{\D}$ et $\R$
  respectivement. 
\end{lemme}

L'intérêt de la fonction $\wdt r$ est d'être
semi-continue supérieurement en les deux variables. Notons que $r\leq \wdt r$
et $R\leq \wdt R$. Il s'agit d'exhiber des valeurs de $\alpha\in \R$ telles que
$r(\alpha,e^{i\theta})=0$ pour toute valeur de $\theta$. Plus précisément, 
\begin{prop}\label{exemples}
Il existe un $G_\delta$ dense $E$ de $\R$ tel que $r(\alpha,e^{i\theta})=0$
pour tout $\alpha\in E$ et tout $\theta \in \R$.
\end{prop}
\noindent \underline{Preuve :} Commençons par montrer que 
\begin{equation}\label{1}
\wdt R=0 \text{ sur } \Q.
\end{equation}
Soit $\alpha=p/q\in \Q$ et $w_0=e^{i2\pi\beta},\; \beta\in \Q$. L'itérée
$q$-ième de $f_\alpha$ est de la forme $\wdt f_\alpha
(w,z):=f_\alpha^q(w,z)=(w,P_w(z))$ où 
$P_w(z)=w^qz+\dots +z^{2^q}$. L'origine étant dans le Julia du polynôme
$P_{w_0}$, elle adhère à son bassin d'attraction de l'infini. Ceci implique que
$\wdt r(\alpha,w_0)=0$ et donc que $\wdt R(\alpha)=0$.

Vérifions maintenant que 
\begin{equation}\label{2}
r(\alpha,\cdot)\equiv R(\alpha) \text{ pour } \alpha\in \R\priv\Q.
\end{equation}
Fixons pour cela $w$ et $w'$ dans $b\D$ et posons $w_n:=we^{i2\pi n\alpha}$
pour $n\in \N$. Il est clair que 
$$
r(\alpha,w)\leq r(\alpha,w_n)\leq r(\alpha,w_{n+1}).
$$
On en déduit, grâce à la semi-continuité supérieure de $r(\alpha,\cdot)$ et en
choisissant une sous-suite $w_{n_k}$ tendant vers $w'$, que $r(\alpha,w')\geq
\overline \lim \, r(\alpha,w_{n_k})\geq r(\alpha,w)$. Ainsi $r(\alpha,.)$ est
constante et coïncide avec $R(\alpha)$.

Nous sommes maintenant en mesure de définir l'ensemble $E$.
Puisque $\wdt R$ est semi-continue supérieurement, $\{\wdt R=0\}$ est un
$G_\delta$ dont la densité résulte de (\ref{1}). Par le théorème de Baire,
$E:=\{\wdt R=0\}\cap \R\priv \Q$ est également un $G_\delta$ dense. D'après
(\ref{2}), on a  $r(\alpha,e^{i\theta})=R(\alpha)\leq \wdt
R(\alpha)=0$ pour tout $\alpha\in E$ et tout $\theta\in \R$.\cqfd

\section{Préliminaires.}
\subsection{Structure des applications holomorphes propres.}
Nous présentons ici les propriétés connues des applications holomorphes propres
qui nous seront utiles. Si $F$ est une application holomorphe
propre entre deux domaines de $\C^k$, nous noterons $V_F$ son lieu de
branchement. C'est l'ensemble d'annulation du Jacobien de $F$.

La première propriété, pour laquelle nous renvoyons à \cite{rudin2},
 décrit la structure de ces applications en terme de revêtements :
\begin{theoreme}\label{revet}
Si $F:\om_1\lra \om_2$ est une application holomorphe propre entre deux domaines
bornés de $\C^k$ alors $F$ est un revêtement ramifié fini, c'est-à-dire :
\begin{enumerate}
\item[i)] $F$ est surjective et ouverte,
\item[ii)] $V_F$ et $F(V_F)$ sont des ensembles analytiques de codimension $1$,
\item[iii)] $F:\om_1\priv F^{-1}(F(V_F))\lra \om_2\priv F(V_F)$ est un
  revêtement fini de degré $d$. Pour tout $z\in \om_2$, Card$\,F^{-1}(z)<d$ si
  et seulement si $z\in F(V_F)$.
\end{enumerate}
\end{theoreme}
Ainsi, lorsque le domaine $\om_2$ est simplement connexe, il suffit pour que
$F$ soit injective qu'elle ne branche pas ($V_F=\emptyset$). Pour les
auto-applications, la simple connexité peut être remplacée par une hypothèse
de régularité du bord :
\begin{lemme}[Pinchuk \cite{pinchuk}]\label{pincuk}
Si $F:\om\lra \om$ est une auto-application holomorphe propre d'un
domaine borné à bord lisse alors $F$ est un  automorphisme de
$\om$ si et seulement si  $V_F=\emptyset$.  
\end{lemme}
Notre approche consistant à étudier le système dynamique engendré par la
restriction de $F$ à $b\om$ ou tout au moins à $\spc(b\om)$, nous utiliserons
les résultats de prolongement suivants :
\begin{theoreme}[Bell \cite{bell}]\label{prolong1}
Toute application  holomorphe propre entre domaines pseudoconvexes bornés de
 $\C^k$ à bords lisses se prolonge différentiablement à $\spc(b\om_1)$. 
\end{theoreme}
\begin{theoreme}[Boas-Straube \cite{boasstraube}]\label{prolong}
Soit $F:\om_1\lra\om_2$ une application holomorphe propre entre domaines
pseudoconvexes bornés de $\C^k$ à bords lisses.
Si $\om_1$ possède une fonction \psh\ lisse et définissante globale alors $F$
s'étend en une application lisse de $\overline{\om_1}$.
\end{theoreme}
Diederich et Fornaess ont montré pourquoi la stricte pseudoconvexité
est une obstruction au branchement (voir \cite{diederichfornaess1} lemme 4). Il
s'agit d'une observation fondamentale que nous incorporons à la proposition
suivante (assertion ii)): 
\begin{prop}\label{strucbord}
Soit $F:\om_1\lra \om_2$ une application holomorphe propre entre domaines
pseudoconvexes bornés de $\C^k$ à bords lisses. Si $F\in {\cal
  C}^\infty(\overline{\om_1})$ alors :
\begin{enumerate}
\item[i)] $F:b\om_1\lra b\om_2$ est surjective, 
\item[ii)] $F(\spc(b\om_1))\subset \spc(b\om_2)$ et $\overline{V_F}\cap
  \spc(b\om_1)=\emptyset$.
\item[iii)] $\eta\in\fpc(b\om_1)$ et $F(\eta)\in \spc(b\om_2) \Longrightarrow
  \eta\in \overline{V_F}$. 
\end{enumerate}
\end{prop} 
Il est facile de préciser la structure des restrictions au bord des
applications holomorphes propres qui fixent le lieu de faible
pseudoconvexité. Ceci nous sera utile dans l'étude des domaines disqués.   
\begin{prop}\label{revetbord}
Soit $F:\om_1\lra \om_2$ vérifiant les hypothèses de la proposition
\ref{strucbord}. Si de surcroît $F(\fpc(b\om_1))\subset \fpc(b\om_2)$ alors   
\begin{enumerate}
\item[i)] $F^{-1}(\fpc(b\om_2))=\fpc(b\om_1)$, $F(\fpc(b\om_1))=\fpc(b\om_2)$
  et
\item[ii)] $F:\spc(b\om_1)\lra \spc(b\om_2)$ est un revêtement fini de même
  degré que $F_{|\om_1}$. 
\end{enumerate}
\end{prop}
Expliquons brièvement le point ii). Notons $d$ le degré de $F$ dans $\om$. 
Puisque l'application $F:\spc(b\om_1)\lra \spc(b\om_2)$ est un difféomorphisme
local, il suffit de vérifier que Card$\,F^{-1}(p)=d$ pour tout $p\in
\spc(b\om_2)$. Montrons que 
\begin{equation}\label{revetbord1}
\text{Card}\, F^{-1}(p)\geq d \hspace{1cm} \forall p\in \spc(b\om_2).
\end{equation} 
Soit $p\in \spc(b\om_2)$. Puisque $F(\fpc(b\om))\subset \fpc(b\om)$, la
deuxième assertion de la proposition précédente montre que $p\notin
\overline{F(V_F)}$. Il existe donc une suite $p_n$  
de $\om_2 \priv F(V_F)$ tendant vers $p$. Chaque $p_n$ a $d$ préimages
distinctes par $F$ que nous noterons 
$x_n^i$, $1\leq i\leq d$. Comme $F$ est continue sur $\om_1$, $(x_n^i)_n$ tend
vers $b\om_1$ pour tout $i$. Toute valeur d'adhérence des suites $(x_n^i)_n$
est une préimage de $p$. Comme $p\notin \overline{F(V_F)}$, $x_{\phi(n)}^i$ et
$x_{\phi(n)}^j$ ne peuvent avoir même limite pour $i\neq j$. Le nombre de
valeurs d'adhérence des $(x_n^i)_n$ est donc au moins $d$, ce qui établit
(\ref{revetbord1}). L'inégalité Card$\,F^{-1}(p)\leq d$ s'obtient facilement
de manière analogue.
L'application $F:\spc(b\om_1)\lra \spc(b\om_2)$ est donc un difféomorphisme
local à fibres finies de cardinal constant égal à $d$. Il s'ensuit que $F$
est un revêtement de degré $d$.

\subsection{Distance $\text{CR}$.}
Nous utiliserons une distance adaptée à la géométrie
$CR$ introduite par Nagel, Stein et Wainger \cite{nsw}. 
Soit $S$ une sous-variété lisse de $\C^k$, $x$ et $y$ deux points de $S$. On
appelle chemin complexe entre $x$ et $y$ tout chemin $\gamma :[0,l]\lra S$
joignant $x$ à $y$, $\mathcal C^1$ par morceaux, et tel que
$\dot{\gamma}(t)\in T^\C_{\gamma(t)} S$ là où cela fait sens. On note
$\ell(\gamma)$ la longueur euclidienne d'un tel chemin. 
Pour $x,y\in S$, on définit la distance $\text{CR}$ entre $x$ et $y$ par :
$$
d^{\text{\tiny CR}}_S(x,y):=\text{inf}\{\ell(\gamma),\;\gamma \text{ chemin
  complexe entre $x$ et $y$}\}. 
$$
En l'absence d'ambiguïté, on ne précisera pas la dépendance en $S$ de
$d^{\text{\tiny CR}}$. Les boules correspondantes de centre $x$ et de
rayon $\delta$ sont notées $B^{\text{\tiny CR}}(x,\delta)$. 

La proposition suivante est bien connue (voir par exemple \cite{nsw}).
\begin{prop}\label{minimale}
Soit $S$ une hypersurface lisse minimale de $\C^k$ (i.e. ne contenant pas de
germe d'hypersurface holomorphe).  La topologie associée à la
distance $\text{CR}$ sur $S$ coïncide avec la toplogie usuelle sur $S$.
En particulier, toute hypersurface compacte connexe et minimale est
$d^{\text{\tiny CR}}$-bornée. 
\end{prop}
\noindent En voici une preuve rapide.
Il est clair que la distance $\text{CR}$ est supérieure à la distance euclidienne, ou
encore que $B^{\text{\tiny CR}}(x,\delta)\subset B(x,\delta)$. Il suffit donc de
trouver $\delta'$ tel que $B^{\text{\tiny CR}}(x,\delta)\supset B(x,\delta')$.
Si $S$ est minimale, la distribution d'hyperplans $T^\C S$ dans $TS$ n'est
intégrable en aucun point. En effet, une hypersurface intégrale de cette
distribution de plans serait une variété de $\C^k$ dont les plans tangents
sont complexes et serait holomorphe d'après un théorème de Levi-Civita
(\cite{chabat}). La boule $B^{\text{\tiny CR}}(x,\delta)$ est alors ouverte
pour  la topologie usuelle sur $S$ et contient donc une boule $B(x,\delta')$
si $\delta'$ est suffisamment petit.\cqfd
Nous utiliserons aussi une version plus sophistiquée de cette proposition
\footnote{Joël Merker m'a informé de l'existence de ce
  résultat et je l'en remercie.}
(voir \cite{bernahu}):
\begin{prop} \label{-min} Toute hypersurface de $\C^k$ compacte connexe sans bord 
  est $d^{\text{\tiny CR}}$-bornée.
\end{prop}

L'intérêt principal de cette distance est d'être dilatée
par toute  application $\text{CR}$ dont l'application tangente est dilatante sur $T^\C
S$. Précisément, 
\begin{prop}\label{dilatc}
Soient $S_1$ et $S_2$ deux hypersurfaces lisses et minimales de $\C^k$. Soit
$F:S_1\lra S_2$ une application $\text{CR}$ qui est un difféomorphisme local en tout
point de $S_1$. Alors, pour toute boule $B^{\text{\tiny CR}}(x,r)$ relativement compacte
dans $S_1$, on a :
$$
\left(\Vert F'(q)u\Vert\geq C\Vert u \Vert, \hspace{0.2cm}\forall q\in
  B^{\text{\tiny CR}}(x,r), \forall   u\in T_q^\C S_1\right)\Longrightarrow 
F\big(B^{\text{\tiny CR}}(x,r)\big)\supset B^{\text{\tiny CR}}(F(x),Cr).   
$$  
\end{prop}
\noindent \underline{Preuve :} Si $F(B^{\text{\tiny CR}}(x,r))=S_2$ la propriété est vraie,
on suppose donc que $F(B^{\text{\tiny CR}}(x,r))\neq S_2$. Comme $F$ est un
difféomorphisme local, $F$ est ouverte et il suffit de voir que  
$$
d^{\text{\tiny CR}}\big(F(x),bF(B^{\text{\tiny CR}}(x,r))\big)\geq Cr.
$$ 
Soit donc $z\in bF(B^{\text{\tiny CR}}(x,r))$ et $\gamma$ un chemin complexe entre $F(x)$
et $z$ paramétré par la longueur d'arc, de longueur $\ell(\gamma)\leq
d^{\text{\tiny CR}}(F(x),z)+\eps$.  
Comme $F$ est un difféomorphisme $\text{CR}$ local en tout point de $S_1$, on
relève $\gamma$ en un chemin complexe  $\wdt \gamma$ tant que $\wdt\gamma$ ne
sort pas de $B^{\text{\tiny CR}}(x,r)$. Autrement dit, il existe $\wdt l\leq
\ell(\gamma)$ et $\wdt \gamma:[0,\wdt l]\lra \overline{B^{\text{\tiny
      CR}}(x,r)}$ joignant $x$ à $bB^{\text{\tiny CR}}(x,r)$  tels que
$F\circ \wdt \gamma(t)=\gamma(t)$ pour tout $t\in[0,\wdt l]$.
Alors 
$$
\begin{array}{rl}
d^{\text{\tiny CR}}(F(x),z)+\eps \geq& \ell(\gamma) \geq \wdt l = \int_0^{\wdt l}\Vert
 \dot{\gamma}(t)\Vert dt
 =\int_0^{\wdt l} \Vert F'(\wdt \gamma(t))\dot{\wdt\gamma}(t)\Vert dt\\
 \geq &  C \int_0^{\wdt l}\Vert \dot{\wdt\gamma}(t)\Vert dt \geq C\ell(\wdt
 \gamma)\geq Cr \text{ puisque } \wdt \gamma(\wdt l)\in bB^{\text{\tiny
 CR}}(x,r).  
\end{array} 
$$
En faisant tendre $\eps$ vers $0$, on obtient bien $d^{\text{\tiny CR}}(F(x),z)\geq
Cr$. \cqfd

Concluons cette partie en précisant l'allure de $B^{\text{\tiny CR}}_{b\om}(x,\eps)$
lorsque $\om$ est un domaine de $\C^k$ à bord lisse et $p$ un point de stricte
pseudoconvexité de son bord. Il s'agit d'un cas particulier très simple du
théorème 4 de \cite{nsw}. 
\begin{prop}\label{boulecplxe}
Soit $p\in \spc(b\om)$. Il existe deux fonctions positives $\tau_1,\tau_2$
définies sur $\R^*_+$ et deux constantes $\tau_0,\eps_0>0$ telles que 
$$
\forall \tau<\tau_0, \; \forall \eps\leq \eps_0,\; U_{\eps,\tau_1(\tau)}\subset
B^{\text{\tiny CR}}(p,\tau\sqrt \eps)\subset U_{\eps,\tau_2(\tau)}.
$$
De plus, on peut prendre $\tau_2$ tendant vers $0$ avec $\tau$.
\end{prop}
Cet énoncé permet de substituer aux boules $\text{CR}$ des objets plus
maniables que nous allons décrire précisément.

Pour $p\in b\om$, $\vec N(p)$ désigne la normale rentrante à $b\om$ en
$p$. Quitte à changer de coordonnées, le domaine $\om$ est strictement
{\textit{convexe}} sur un voisinage $O$ de $p$, et est donné par  
$$
\begin{array}{l}
\om\cap O=\{(w,z)\in \C\times \C^{k-1}\textit{ t.q. }\re w\geq |z|^2+a(\im w,z)\}\\
\text{ où } a(v,z)\in O(v^2,|z|^3)
\end{array}
$$ 
(dans ces coordonnées, si $w=u+iv$, on a $p=(0,0)$ et $\vec
N(p)=\displaystyle \partial/\partial u$).\\
Comme $\om\cap O$ est convexe, on peut définir une projection $\pi$ de $O\cap
\om$ sur $O\cap b\om$ dans la direction $\vec N(p)$. Elle est donnée dans les
coordonnées ci-dessus par : 
$$
\pi(u+iv,z)=(|z|^2+a(v,z)+iv,z).
$$
Rappelons que $T_p b\om=T_p^\C b\om\oplus$Vect$_\R (i\vec N(p))$. 
Nous notons $B^{t,\C}(p,r)$ la boule de centre $0$ et de rayon $r$ dans
$T^\C_p b\om$ et $B^{t,\R}(p,r)$ le segment de centre $0$ et de rayon $r$ dans
Vect$_\R(i\vec{N}(p))$. Ceci étant posé nous sommes en mesure de formuler la  
\begin{definition}\label{defueps}
$$
\begin{array}{l}
\ds p_\eps:=p+\eps\vec N(p)\\
\ds F_{\eps,\tau}:=p_\eps+B^{t,\C}(p,\tau \sqrt\eps)\times B^{t,\R}(p,\tau
\eps)\\
\ds U_{\eps,\tau}:=\pi(F_{\eps,\tau}).
\end{array}
$$
\end{definition}
Signalons dès à présent que $\eps$ est destiné à tendre vers $0$ alors que $\tau$
est un paramètre de contrôle qui sera fixé. Ces ensembles sont définis
dès que $\eps$ est suffisamment petit et $\tau$ plus petit que $1$.
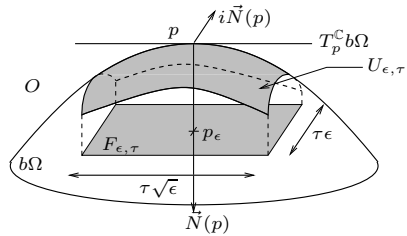
\begin{figure}[h]
\begin{center}
\input uepsilon3.pstex_t
\end{center}
\end{figure}


\subsection{Quelques estimations de la métrique de Kobayashi.}
Nous utiliserons les estimations usuelles suivantes :
\begin{lemme}\label{kobaya}
Soit $\om$ un domaine borné à bord lisse de $\C^k$ et $p$ un point de stricte
pseudoconvexité de $b\om$. On reprend les coordonnées et les notations
introduites à la fin de la partie précédente.

Il existe $\tau_0,\eps_0>0$ tels que pour tout $\tau<\tau_0$ et $\eps<\eps_0$ :
\begin{enumerate}
\item[i)] pour $z\in F_{\eps,\tau}$, le disque $\Delta(z)$ de centre $z$ et
  de rayon $\Vert z-\pi(z)\Vert$ dans $z+$Vect$_\C(z-\pi(z))$ est inclus dans $\om$,
\item[ii)] $Diam_{K_\om}(F_{\eps,\tau})\leq \psi_1(\tau)$ où
  $\psi_1(\tau)$ tend vers $0$ avec $\tau$,
\item[iii)] pour $z\in F_{\eps,\tau}$, Diam$_{K_\om}\big(\mathcal
  C^\alpha(z,\pi(z))\big)\leq \psi_2(\alpha)$ où $\psi_2(\alpha)$ tend vers $0$
  avec $\alpha$. 
\end{enumerate} 
\end{lemme} 
\noindent \underline{Preuve : \\ }
\begin{enumerate}
\item[i)] C'est évident car le bord de $\om$ est lisse.
\item[ii)] Effectuons une remise à l'echelle en $p$ de paramètre $\eps$, {\it i.e.}
considérons le domaine $\om_\eps=\Lambda_\eps(\om)$ où
$\Lambda_\eps(w,z)=(\frac{w}{\eps},\frac{z}{\sqrt\eps})$. On a  
$$
F_{1,\tau}:=\Lambda_\eps(F_{\eps,\tau})=\{(w,z)\in\C^k,\; 
\re w =1,\; |\im w|\leq \tau, |z|\leq\tau \}.
$$
Un calcul très simple (voir \cite{pinchuk2}) montre que $\wdt \om_\eps$
tend vers $\Sigma:=\{\re w\geq |z|^2\}$ au sens de Hausdorff lorsque
$\eps$ tend vers $0$. Comme $B:=B\big((1,0),1/2\big)\Subset\Sigma$, il
s'ensuit que $B\subset\wdt \om_\eps$ pour $\eps\leq \eps_1$ lorsque $\eps_1$
est suffisamment petit. Puis $F_{1,\tau}\subset B$ pour $\tau\leq\tau_2$ où
$\tau_2$ est assez petit. Ainsi, pour $\tau<\tau_2$ et $\eps<\eps_1$, on a 
$$
\text{Diam}_{K_\om}(F_{\eps,\tau})=
\text{Diam}_{K_{\om_\eps}}(F_{1,\tau})\leq
\text{Diam}_{K_B}(F_{1,\tau})
$$ 
et il suffit de prendre $\psi_1(\tau):=\,$Diam$_{K_B}(F_{1,\tau})$.
\item[iii)] D'après i), $\Delta(z)$ est inclu dans $\om$ pour $\eps\leq
  \eps_1$. Donc  
$$
\text{Diam}_{K_\om}\big(\mathcal C^\alpha(z,\pi(z))\big) \leq
\text{Diam}_{K_{\Delta(z)}}\big(\mathcal C^\alpha(z,\pi(z))\big)
\leq \text{Diam}_{K_{\D_1}}(\D_\alpha)
$$
et il suffit de prendre $\psi_2(\alpha):=\,$Diam$_{K_{\D_1}}(\D_\alpha)$.
\end{enumerate}\cqfd
\subsection{Lemme de Hopf.}
Le très classique lemme de Hopf nous permettra d'estimer les dérivées normales
des applications holomorphes propres. En voici une version adaptée à nos
besoins 
: 
\begin{lemme}[Hopf]\label{hopf}
Soit $\chi\in \mathcal C^{\infty}(\overline{\D_r})$ une fonction
sous-harmonique et négative sur $\D_r$ telle que $\chi(-r)=0$ et 
$\chi(re^{i\theta})\leq -1$ pour $\theta\in [-\alpha,\alpha]$. 
Alors $\ds\chi(-r+t)\leq -\frac{\alpha t}{4\pi r}$ pour $t\leq r$ et donc $\ds
\Vert \grad\chi(-r)\Vert\geq \frac{\alpha}{4\pi r}$.
\end{lemme}
\noindent \underline{Preuve :} Il s'agit d'un calcul facile sur le noyau de Poisson.
$$
\begin{array}{rl} 
\chi(-r+t)  \leq &\ds\frac{1}{2\pi}\int_{-\pi}^\pi
\frac{r^2-(r-t)^2}{|-r+t-re^{i\theta}|^2} \chi(re^{i\theta})d\theta \\
  \leq &  \ds -\frac{1}{2\pi}\int_{-\alpha}^\alpha
\frac{r^2-(r-t)^2}{|-r+t-re^{i\theta}|^2}\, d\theta 
 \leq \ds -\frac{\alpha}{\pi}\frac{t(2r-t)}{4r^2} 
 \leq  \ds -\frac{\alpha t}{4\pi r}.
\end{array}
$$
\cqfd
\section{Une situation modèle : le bassin d'un point attractif.}
Il s'agit du cas où la dynamique de l'application est récurrente et le rétract
de dimension nulle. Ceci est impossible lorsque le bord du domaine est trop
régulier. Plus précisément,
\begin{theoreme}
Un domaine borné de $\C^k$ à bord lisse minimal
possédant une fonction \psh\ définissante globale $\chi$ et dont le bord est minimal
ne peut être le bassin d'attraction de l'un de ses points pour une
transformation holomorphe.
\end{theoreme}
La démonstration permet d'illustrer les méthodes utilisées dans cet article
tout en évitant les difficultés techniques liées au cas général. En voici le
principe. Soit $f:\om\lra \om$ une application holomorphe propre telle que
$(f^n)_n$ converge localement uniformément vers $a\in \om$. Il suffit pour
aboutir à une contradiction de montrer que $f$ est un automorphisme
c'est-à-dire que $f$ est injective. On voit dans un premier temps que la suite
$f^n$ présente un "défaut" d'équicontinuité dans la direction normale au
voisinage de $b\om$.
On en déduit que les dérivées de $f^n$ évaluées sur les directions
tangentielles  complexes explosent sur $\spc(b\om)$. Il s'ensuit que 
$b\om=\spc(b\om)$ et donc que $f$ est injective. Voyons cela plus en détail.

Introduisons quelques notations. Soit $U$ un ouvert de $b\om$ relativement
compact dans $\spc(b\om)$, $p$ un point de $U$ et $\tau>0$ un réel tel que
$B^{\text{\tiny CR}}(p,\tau)\subset U$. Rappelons que $f$ se prolonge en une
application lisse de $\overline{\om}$ (théorème \ref{prolong}). Pour $q\in
b\om$, nous notons
$$
n_q(f^n):=\langle {f^n}'(q)\vec N(q),\, \vec N(f^n(q))\rangle.
$$ 
Cette quantité mesure le "taux d'échappement dans la direction normale" de
$f^n$.\vspace{0.2cm}\\ 
\textbf{Fait 1 :}\textit{ La suite $n_q(f^n)$ tend vers l'infini
uniformément sur $b\om$.}\vspace{0.2cm}\\
\underline{Preuve :} Soit une coquille $A_\eps:=\bigcup_{q\in
  b\om}B(q_\eps,\eps/2)$ d'épaisseur $\eps/2$. Pour $n\geq n_\eps\gg 1$, on a
$f^n(A_\eps)\subset \{\chi\leq \chi(a)/2\}$. Le lemme de Hopf
appliqué à $\chi\circ f^n_{|\D(q,q_\eps)}$ montre que $\Vert\grad \chi\circ f^n\Vert\gtrsim
1/\eps$ pour $n\geq n_1$. Ceci suffit car
$n_q(f^n)\simeq  \Vert\grad \chi\circ f^n\Vert$ :
$$
\Vert\grad \chi\circ f^n\Vert=\langle \grad \chi\circ f^n(q),\vec N(q)\rangle=\langle \grad
\chi(f^n(q)),{f^n}'(q)\vec N(q)\rangle \simeq n_q(f^n).
$$
\cqfd
\textbf{ Fait 2 :}\textit{
Il existe une constante $c>0$ telle que l'inégalité} 
$$
\Vert F'(q)u\Vert\geq cn_q(F)^{1/2}\Vert u\Vert \hspace{0.5cm} \forall q\in
U,\forall u\in T^\C_q b\om
$$
\textit{ait lieu pour toute auto-application $F$ holomorphe propre de
  $\om$.}\vspace{0.2cm}\\
\underline{Preuve :} 
Comme $U\Subset \spc(b\om)$, il existe une constante $\kappa$ telle que pour
toute fonction \psh\ définissante globale $\psi$ de $\om$,
$$
\L(\psi,q,u)\geq \kappa \Vert\grad \psi(q)\Vert\, \Vert u\Vert^2 \hspace{0.5cm}
\forall q\in U,\; \forall u\in T^\C_q b\om.
$$  
Nous avons remarqué que $\Vert\grad \chi\circ F\Vert\simeq n_q(F)$ lorsque $F$
est une auto-application holomorphe propre de $\om$. Ainsi, $\L(\chi\circ
F,q,u)\gtrsim n_q(F)\Vert u\Vert^2$. Or $\L(\chi\circ
F,q,u)=\L(\chi,F(q),F'(q)u)\lesssim \Vert F'(q)u\Vert^2$ puisque $b\om$ est
lisse. On obtient donc finalement :
$$
\Vert F'(q)u \Vert\gtrsim n_q(F)^{1/2}\Vert u\Vert.
$$\cqfd
\textbf{Conclusion :} On voit à l'aide des faits 1 et 2 que $\Vert
{f^n}'(q)u\Vert$ tend vers l'infini uniformément sur 
$\{(q,u)\in U\times T^\C_q b\om,\; \Vert u\Vert=1\}$. La proposition
\ref{dilatc} montre alors que $f^n(B^{\text{\tiny CR}}(p,\tau))$ contient 
des boules $\text{CR}$ aussi grandes que souhaité pourvu que $n$ soit
assez grand. Comme $b\om$ est $d^{\text{\tiny CR}}$-bornée (proposition \ref{-min}),
il existe un entier $n_0$ tel que $f^{n_0}(B^{\text{\tiny CR}}(p,\tau))
=b\om$. D'après la propriété \ref{strucbord}.ii), $b\om$ est strictement 
pseudoconvexe et $f$ est un biholomorphisme. \cqfd

\section{Dérivées au bord des applications holomorphes propres.}
L'objet de cette partie est d'estimer les dérivées au bord d'auto-applications
holomorphes propres à partir de leur dynamique à l'intérieur du domaine.
\subsection{Dérivées normales.}
Nous estimons ici les dérivées normales. 
Plus précisément, nous minorons les quantités  
$$
n_q(f^n):=\langle {f^n}'(q)\vec N(q),\vec N(f^n(q))\rangle
$$
en fonction de la dynamique de $f$.
Le lemme de Hopf fournit de telles estimations 
lorsqu'il existe des "barrières" \psh\ appropriées sur le
domaine. 

Dans l'étude de la dynamique récurrente (cinquième partie), nous supposerons
l'existence de fonctions définissantes globales \psh. Nous procéderons comme
dans la preuve du Fait 1 de la situation modèle. Cependant, lorsque le rétract
est de dimension positive, l'existence d'un compact contenant
les images $f^n(A_\eps)$ des coquilles $A_\eps$ pour tout $\eps>0$ et tout
$n\geq n(\eps)\gg 1$ n'est plus assurée. Tout au plus trouverons nous une
suite de points $p_\eps$ tendant normalement vers $p\in \spc(b\om)$ et telle que
$f^{n}(p_\eps)\in K\Subset \om$ pour $n\geq n(\eps)$. Dans ces conditions,
nous ne contrôlerons plus $n_q(f^n)$ pour $q\in b\om$ mais
seulement sur des boules $\text{CR}$ du type $B^{\text{\tiny CR}}(p,\tau\sqrt
\eps)$. Ceci fait l'objet de la proposition suivante.
\begin{prop}\label{norm1}
Soit $\om\subset \C^k$ un domaine possédant une fonction \psh\ lisse et
définissante globale  $\chi$. Soit $p\in \spc(b\om)$. Soit $F:\om\lra \om$ une 
application holomorphe propre lisse sur $\overline{\om}$. 
Il existe trois constantes $\eps_0,c,\tau>0$ (indépendantes de $F$) telles que
pour $\eps\leq\eps_0$ :
$$
\chi\big(F(p_\eps)\big)\leq -1\Longrightarrow n_q(F):=\langle F'(q)\vec
N(q),\vec N(F(q))\rangle\geq\frac{c}{\eps} \hspace{0.3cm}\forall q\in
B^{\text{\tiny CR}}(p,\tau\sqrt\eps).
$$
\end{prop}
\noindent \underline{Preuve :} On peut supposer $b\om$ convexe en $p$ et 
utiliser les notations introduites dans la definition \ref{defueps}.
Soit $\eps_0,\tau_0,\psi_1,\psi_2$ donnés par
le lemme \ref{kobaya} et 
$\eps\leq \eps_0$. Supposons $\chi\left[F(p_\eps)\right]\leq -1$. Soit
$\delta>0$ tel que  
$$
\chi(z)\leq -1 \text{ et } d_K(x,z)\leq \delta\Longrightarrow \chi(x)\leq
-1/2.
$$
Soient $0<\tau_2\leq \tau_0$ et $\alpha>0$ tels que
$\psi_1(\tau_2)+\psi_2(\alpha)<\delta$. Fixons $z\in F_{\eps,\tau_2}$.
D'après le lemme \ref{kobaya}, la distance de Kobayashi entre $p_\eps$ et tout
point de $C^\alpha(z,\pi(z))$  est inférieure à $\delta$. Puisque $F$
contracte cette distance :
$$
d_K\big(F(p_\eps),z'\big) < \delta \hspace{0.5cm}\forall z'\in F\big(\mathcal
C^\alpha(z,\pi(z))\big) 
$$
et 
$$
\chi\circ F\leq -\frac{1}{2} \text{ sur }
\mathcal C^\alpha(z,\pi(z)).
$$
Toujours d'après le lemme \ref{kobaya}, $\D(z,\pi(z))\subset \om$. Comme
$F$ est holomorphe propre,  
$\chi\circ F$ est une fonction \psh\  définissante de $\om$. Appliquons alors
le lemme de Hopf \ref{hopf} à $\chi\circ  F_{|\D(z,\pi(z))}$. Pour cela,
notons  $q:=\pi(z)$, $r:=\frac{\Vert z-q\Vert}{2}$ et $v:=\frac{z-q}{\Vert
  z-q\Vert}$. Observons que $q\in U_{\eps,\tau_2}$. On a 
$$
\frac{\alpha}{4\pi r}\leq |d[\chi\circ F]_q\cdot
v|=|d\chi(F(q))\cdot F'(q)v|=|\,\langle\grad\chi(F(q)),F'(q)v \rangle\,|.
$$  
En remarquant que $r\leq \eps/2$ (par convexité) et en écrivant
$v=a\vec N(q)+u$ où $0\leq a\leq 1$ et $u\in T_q b\om$, il vient :
$$ 
\begin{array}{ll}
\left|\langle\grad \chi(F(q)),F'(q)v\rangle\right| &=\Vert \grad
 \chi(F(q))\Vert \,\langle\vec  N(F(q))\,,\,a F'(q)\vec N(q)+F'(q)u\rangle\\
 &= a\Vert \grad \chi(F(q))\Vert \,\langle\vec N(F(q))\, ,\, F'(q)\vec N(q)\rangle .
\end{array}
$$

On a obtenu :
$$
n_q(F)\geq \frac{\alpha}{2\pi\eps\Vert
  \grad\chi\Vert_\infty}=\frac{c}{\eps}\, ,\hspace{0.5cm} \forall q\in U_{\eps,\tau_2}.
$$
La constante $c$ est indépendante de $F$. Il suffit de prendre $\tau$ tel que
$\tau_2(\tau)\leq \tau_2$ pour avoir le résultat annoncé puisqu'alors
$B^{\text{\tiny CR}}(p,\tau\sqrt\eps)\subset U_{\eps,\tau_2}$ en vertu de la
proposition \ref{boulecplxe}.    \vspace{0.2cm}\cqfd 
\indent Dans l'étude de la dynamique non-récurrente (sixième partie), nous
utiliserons des fonctions pics \psh\  et la 
\begin{prop}\label{norm2}
Soit $\om$ un domaine pseudoconvexe à bord lisse de $\C^k$ et $F:\om\lra
\om$ une application holomorphe propre lisse sur $\overline{\om}$. Soit
$p\in b\om$. Supposons qu'il existe une fonction $\chi\in PSH(\om)\cap \mathcal
C^0(\overline{\om})$ telle que 
\begin{description}
\item{\sbull} $\chi\leq 0$,
\item{\sbull} $F\big(\mathcal C_\alpha(p_\eps,p)\big)\subset \{\chi\leq
  -1\}$,
\item{\sbull} $\chi$ est nulle sur un voisinage de $F(p)$ dans $b\om$,
\item{\sbull} $\chi$ est $L$-lipschitzienne en tout point d'un voisinage de
  $F(p)$ dans $b\om$.  
\end{description}
Alors $\ds n_p(F)\geq \frac{\alpha}{4\pi L\eps}$. 
\end{prop} 
\noindent\underline{Preuve :} 
On peut prolonger $\chi$ à un voisinage de $b\om$ tout en la conservant 
$L$-Lipschitzienne sur un voisinage de $F(p)$ dans $b\om$.
Le lemme de Hopf \ref{hopf} appliqué à $\chi\circ F$ dans $\D(p_\eps,p)$
montre que : 
$$
\chi\circ F(p+t\vec N(p))\leq -\frac{\alpha}{4\pi\eps}\,t. 
$$
Par ailleurs, 
$$
\begin{array}{ll}
F(p+t\vec N(p))&=F(p)+tF'(p)\vec N(p)+O(t^2)\\
 & =F(p)+t n_p(F)\vec N(F(p))+t\vec R+O(t^2)
\end{array}
$$ 
où $\vec R$ est la projection orthogonale de $F'(p)\vec N(p)$ sur
$T_{F(p)}b\om$. Il existe clairement un point $q_t$ de $b\om$ tel que
$F(p)+t\vec R=q_t+O(t^2)$.  
Si $t$ est suffisamment petit, $\chi$ est $L$-lipschitzienne en $q_t$. On a donc : 
$$
\begin{array}{ll}
\ds\frac{\alpha}{4\pi\eps}t \leq |\chi\circ F(p+t\vec
N(p))|&=|\chi\big(q_t+t n_p(F)\vec N(F(p))+O(t^2)\big)|\\
 & \leq |\chi(q_t)|+ Lt n_p(F)+O(t^2) =Lt n_p(F)+O(t^2).
\end{array}
$$  
On conclut en faisant tendre $t$ vers $0$.\cqfd
\subsection{Dérivées tangentielles complexes.}
Il s'agit de transférer les estimations des dérivées normales aux dérivées
tangentielles complexes. On utilise pour cela un argument standard reposant
sur la stricte pseudoconvexité et la fonctorialité de la forme de Levi. Dans
la situation modèle, ceci correspondait au Fait 2. 
\begin{lemme}\label{normcomp}
Soient $\om_1$ et $\om_2$ deux domaines pseudoconvexes bornés de $\C^k$ à
bords lisses et $K$ un compact de $\spc(b\om_1)$. Soit  $F:\om_1\lra \om_2$
une application holomorphe propre lisse sur $\overline{\om_1}$.
Il existe une constante $c>0$ (ne dépendant que de $K$) telle que :
$$
\Vert F'(q)u\Vert\geq c n_q(F)^{\frac{1}{2}} \Vert u\Vert\, ,  \hspace{0.5cm}\forall q\in
 K \text{ et } \forall u\in T^\C_q b\om.
$$
\end{lemme}
De ce lemme, de la proposition \ref{norm1} et de la proposition \ref{dilatc},
on déduit immédiatement la proposition suivante :
\begin{prop}\label{tanc1}
Soient $\om\subset \C^k$ un domaine borné possédant une fonction  \psh\ lisse
et définissante globale $\chi$. Soit  $p\in 
\spc(b\om)$. Soit $F:\om\lra \om$ une application holomorphe propre lisse sur
$\overline{\om}$. \\
Il existe trois constantes $\eps_0,c,\tau>0$ (indépendantes de $F$) telles que
pour $\eps<\eps_0$ et $\chi~\circ~ F(p_\eps)\leq -1$ :
\begin{enumerate}
\item[i)] $\ds\Vert F'(q)u\Vert\geq \frac{c}{\sqrt \eps}\Vert u\Vert
  \hspace{0.5cm} \forall q\in B^{\text{\tiny CR}}(p,\tau\sqrt \eps),\; \forall u\in T_q^\C b\om$
\item[ii)] $F\big(B^{\text{\tiny CR}}(p,\tau\sqrt \eps)\big)\supset
    B^{\text{\tiny CR}}(F(p),c\tau).$
\end{enumerate}
\end{prop} 
\noindent\underline{Preuve du lemme \ref{normcomp} :}
Soit $c_2$ le maximum des valeurs propres de la forme de Levi de
$b\om_2$  et $c_1$ le minimum des valeurs propres de la forme de Levi de
$b\om_1$ sur $K$. Si $q\in K$ alors $q\in \spc(b\om_1)$ et $F(q)\in
\spc(b\om_2)$. En particulier, on peut trouver deux
fonctions $\chi_1$ et  $\chi_2$ définissant respectivement $b\om_1$ 
et $b\om_2$ au voisinage de $q$ et $F(q)$ telles que $\grad
\chi_1(q)=-\vec N(q)$ et $\grad\chi_2(F(q))=-\vec N(F(q))$.
Comme $F$ est propre, $\chi_2\circ F$ est une fonction définissante de
$b\om_1$ au voisinage de $q$, il existe donc une fonction
$\gamma$ positive et lisse au voisinage de $q$ telle que
$\chi_2\circ F=\gamma \chi_1$. Remarquons que $\gamma(q)=n_q(F)$. En effet, 
$$
\grad[\chi_2\circ F] (q)=\chi_1(q)\grad \gamma (q)+\gamma(q)\grad
\chi_1(q)=-\gamma(q)\vec N(q)
$$
donc 
$$
\begin{array}{ll}
\gamma(q) & = -\langle\grad [\chi_2\circ F](q),\vec N(q)\rangle\\
 & = -\langle\grad \chi_2(F(q)),F'(q)\vec N(q)\rangle\\
 & = \langle\vec N(F(q)),F'(q)\vec N(q))\rangle=n_q(F).
\end{array}
$$
Notons $\L(\phi,a,u)$ la forme de Levi en $a$ d'une fonction $\phi$,
appliquée à un vecteur $u$. Si $u\in T^\C_q b\om_1$,  on a 
$$
\begin{array}{ll}
\L(\chi_2\circ F,q,u)&=\L(\chi_2,F(q),F'(q)u)\\
 & =\L(\gamma \chi_1,q,u).\\
\end{array}
$$
Puisque $\L(\gamma \chi_1,q,u)=\gamma(q)\L(\chi_1,q,u)$
pour tout $q$ de $b\om_1$ et tout $u\in T^\C_qb\om_1$, il vient 
$$
c_2\Vert F'(q)u\Vert^2\geq \L(\chi_2,F(q),F'(q)u)=
\gamma(q)\L(\chi_1,q,u)\geq c_1 n_q(F)\Vert u\Vert^2.
$$
On a donc $\Vert F'(q)u\Vert\geq c n_q(F)^{\frac{1}{2}}\Vert u\Vert$, où
$c=\sqrt{c_1c_2^{-1}}$. \cqfd 
\section{Dynamique récurrente et expansivité au bord.}\label{partierecurrent}
L'objet de cette partie est de démontrer le théorème 1.
\subsection{Lemme des matriochkas.}
Le lemme suivant montre en quoi une dynamique récurrente dans
$\om$ s'accompagne d'un comportement expansif sur les régions strictement
pseudoconvexes de $b\om$. Concrètement, nous ne ferons qu'appliquer la
proposition \ref{tanc1} à une "matriochka" de boules $\text{CR}$ centrées en
$p$.  
\begin{lemme}\label{gale}
Soit $\om\subset \C^k$ un domaine borné possédant
une fonction  \psh\ lisse et définissante globale $\chi$. Soit $p\in
\spc(b\om)$. Soit $f_n:\om\lra \om$ une suite d'applications holomorphes
propres lisses sur $\overline{\om}$.  Supposons qu'il existe des suites  
$t_k\lra 0$ et $n_k\in \N$ telles que $f_n(p_{t_k})\in \{\chi\leq -1\}$ pour
tout $n\geq n_k$. Alors, pour tout $R>0$ et tout voisinage $U$ de $p$, il
existe $n\in \N$ tel que $f_n(U)\supset B^{\text{\tiny CR}}(f_n(p),R)$.
\end{lemme}
\noindent \underline{Preuve :}
Soient $\eps_0,\tau,c$ les constantes données par la
proposition \ref{tanc1}. Fixons un entier $N$ supérieur à $2R/\tau
c$. Choisissons alors successivement 
$\eps_N>\eps_{N-1}>\dots>\eps_1$ parmi les termes de la suite $(t_k)_k$ et  
$n_0\in \N$ de façon à ce que : 
\begin{description}
\item[\sbull] $\B^{\text{\tiny CR}}(p,\tau\sqrt{\eps_N})\subset U$ et $\eps_N<\eps_0$
\item[\sbull] $\ds \eps_i<\frac{\eps_{i+1}}{4}$ 
\item[\sbull] $ f_n(p_{\eps_i})\in \{\chi\leq -1\} \hspace{0.3cm}\forall n\geq n_0$.  
\end{description} 
Pour $n\geq n_0$ donné, posons $F:=f_n$.
Il suffit d'établir par récurrence sur $i$ que :
\begin{equation} \tag{\text{$\mathcal P_i$}} 
F(B^{\text{\tiny CR}}(p,\tau\sqrt{\eps_i}))\supset B^{\text{\tiny
  CR}}(F(p),\frac{ic\tau}{2}) \text{ pour $1\leq i\leq N$.}  
\end{equation}

D'après la proposition \ref{tanc1},
si $q\in B^{\text{\tiny CR}}(p,\tau\sqrt{\eps_1})$ et $u\in T_q^\C b\om$ alors $\Vert
F'(q)u\Vert\geq \frac{c}{\sqrt{\eps_1}}\Vert u\Vert$. En appliquant la
proposition \ref{dilatc}, on a donc :
\begin{equation}\tag{$\mathcal P_1$}
F(B^{\text{\tiny CR}}(p,\tau\sqrt{\eps_1}))\supset B^{\text{\tiny
    CR}}(F(p),c\tau).
\end{equation}

Supposons l'inclusion $(\mathcal P_i)$ satisfaite et établissons
  $(\mathcal P_{i+1})$. Comme $\eps_i<\frac{\eps_{i+1}}{4}$, on a
$B^{\text{\tiny CR}}(x,\frac{\tau\sqrt{\eps_{i+1}}}{2})\subset B^{\text{\tiny
  CR}}(p,\tau\sqrt{\eps_{i+1}})$ 
si $x\in B^{\text{\tiny CR}}(p,\tau\sqrt{\eps_i})$. Il vient donc 
$$
\Vert F'(q)u\Vert\geq \frac{c}{\sqrt{\eps_{i+1}}}\Vert
u\Vert\, ,\hspace{0.5cm}\forall q\in \bigcup_{x\in
  B^{\text{\tiny CR}}(p,\tau\sqrt{\eps_i})}B^{\text{\tiny
    CR}}\big(x,\frac{\tau\sqrt{\eps_{i+1}}}{2}\big),\; \forall u\in  T_q^\C b\om 
$$
et
$$
\begin{array}{c}
\begin{array}{rcl}
F(B^{\text{\tiny CR}}(p,\tau\sqrt{\eps_{i+1}}))\supset & \ds\bigcup_{x\in
  B^{\text{\tiny CR}}(p,\tau\sqrt{\eps_i})} & F\big(B^{\text{\tiny
  CR}}(x,\frac{\tau\sqrt{\eps_{i+1}}}{2})\big) \\ 
 \supset &\ds\bigcup_{x\in
  B^{\text{\tiny CR}}(p,\tau\sqrt{\eps_i})} & B^{\text{\tiny CR}}(F(x),\frac{c\tau}{2}) \\
 \supset & \ds\bigcup_{y\in
 F(B^{\text{\tiny CR}}(p,\tau\sqrt{\eps_i}))} & B^{\text{\tiny CR}}(y,\frac{c\tau}{2}) \\
 \supset & \ds\bigcup_{y\in
  B^{\text{\tiny CR}}(F(p),\frac{ic\tau}{2})} & B^{\text{\tiny CR}}(y,\frac{c\tau}{2}) \\
\end{array} \\
 \hspace{0.9cm} \supset \ds B^{\text{\tiny CR}}\big(F(p),(i+1)\frac{c\tau}{2}\big) .
\end{array}
$$\cqfd

Nous observons maintenant que l'expansivité de la suite
$(f_n)_n$ sur $U$ force la stricte pseudoconvexité de $b\om$ et donc
l'injectivité des $f_n$. 
\begin{prop}\label{gale2}
Lorsque $b\om$ est minimal alors, sous les hypothèses du lemme \ref{gale}, 
les $f_n$ sont des automorphismes de $\om$.
\end{prop}
\noindent \underline{Preuve :} D'après la proposition \ref{-min}, $B^{\text{\tiny
CR}}(\eta,R)$ contient $b\om$ pour tout $\eta$ de $b\om$ pourvu que $R$ soit
supérieur au diamètre $\text{CR}$ de $b\om$. Choisissons un voisinage $U$ de
$p$ assez petit pour que $U\Subset \spc(b\om)$. D'après la proposition
précédente, $b\om\subset 
f_n(U)$ pour $n$ assez grand et donc $b\om$ est strictement pseudoconvexe. Le
lemme \ref{pincuk} montre que les $f_n$ sont des automorphismes. \cqfd 


\subsection{Preuve du théorème 1.}\label{petitedim}
Commençons par rappeler qu'en vertu du théorème \ref{prolong}, $f$ se prolonge
en une application lisse à $\overline{\om}$. Raisonnons par l'absurde. 

Si le rétract est un point, il s'agit d'un point fixe de $f$
  vers lequel $(f^n)_n$ converge localement uniformément. Soit $p$ un point de
  stricte pseudoconvexité de $b\om$ et $K$ un voisinage compact de $a$ dans
  $\om$. Pour tout $\eps$, il existe $n_\eps\in \N$ tel que
  $f^n(p_\eps)\in K$ pour $n\geq n_\eps$.

Lorsque  le rétract est une surface de Riemann $M$, son genre est fini
puisqu'elle est obtenue par une rétraction $\rho$ de $\om$ sur $M$.
Considérons une sous-suite $f^{n_i}$ tendant vers $\rho$. Pour $q\in b\om$,
nous posons 
$$
\rho^*(q):=\bigcap_{t>0}\overline{\{\rho(q+u\vec N(q)),u\leq t\}}.
$$ 
Il s'agit de l'ensemble des valeurs d'adhérences normales de $\rho$ en
$q$. Dans un mémoire sur les fonctions intérieures de la boule, Rudin montre
que les limites radiales de ces fonctions sont denses dans le disque
\cite{rudin}. En adaptant la preuve de ce résultat (voir appendice), nous
établissons le lemme suivant : 
\begin{lemme}\label{limnorm}
Soit $\om$ un domaine pseudoconvexe borné à bord lisse de $\C^k$, $k>1$, et
$M$ une surface de Riemann de genre fini.
Soit $\rho:\om\lra M$ une application holomorphe.
Il existe $p\in \spc(b\om)$ et $x\in M$ tels que $x\in \rho^*(p)$.
\end{lemme} 
\noindent Soit alors $x\in M\subset \om$ le point donné par ce lemme, et
$p_{\eps_k}$ une suite de points telle que $\rho(p_{\eps_k})\lra x$. Comme $f$
a une dynamique récurrente, la suite $(f^n(x))_n$ est relativement compacte dans
$\om$ et il existe $d>0$ tel que $K:=\overline{\cup_{n\geq 0}
  B_{K}(f^n(x),d)}\subset \om$.
Comme $f^{n_i}$ converge vers $\rho$, quitte à extraire, on a 
$f^{n_k}(p_{\eps_k})\in B_K(x,d)$ et $f^n(p_{\eps_k})\in K$ pour $n\geq n_k$.

Dans tous les cas, il existe un point de stricte pseudoconvexité $p$ de
$b\om$, un compact $K$ de $\om$ et des suites $\eps_k\lra 0$ et $n_k\in \N$
telles que $f^n(p_{\eps_k})\in K$ pour $n\geq n_k$. On peut supposer que la
fonction \psh\ définissante $\chi$ est inférieure à $-1$ sur $K$.  
La proposition \ref{gale2} s'applique donc et montre que $f$ est un
biholomorphisme. Ceci est impossible (voir la remarque \ref{pasaut}).\cqfd

\section{Dynamique non-récurrente : transmission de la dynamique à
  $\spc(b\om)$.}
Nous montrons dans cette partie que la non-récurrence de la dynamique de $f$
se transmet à celle de son extension à $\spc(b\om)$. Nous reprenons pour cela
la technique utilisée dans le cas récurrent. Il nous faut donc, en
particulier, construire des fonctions \psh\ négatives sur $\om$, strictement
négatives sur un voisinage d'un point limite de $(f^n)_n$ dans
$\overline{\om}$ et nulles sur de larges portions de $b\om$. Ceci n'est pas
toujours possible. La classe de domaines permettant cette construction a été
caractérisée dans un mémoire de Sibony \cite{sibony}. Ces domaines sont connus
sous le nom de domaines $B$-réguliers. 
\subsection{Domaines $LB$-réguliers.}
\begin{definition}\label{breg}
Un domaine pseudoconvexe borné $\om$ de $\C^k$ à bord lisse est $B$-régulier
si pour tout $p\in b\om$, il existe une fonction $\psi_p\in PSH(\om)\cap
\mathcal C^0(\overline{\om})$ vérifiant $\psi_p(p)=1$ et $\psi_p(z)<1$ pour
$z\neq p$, $z\in \overline{\om}$. 
\end{definition}
Sibony a montré dans \cite{sibony} qu'il existe des fonctions \psh\ de trace au bord 
prescrite sur tout domaine $B$-régulier. En particulier, ces domaines
admettent des fonctions antipics en tout point du 
bord. Elles sont exploitables dans notre cadre de travail
lorsqu'elles sont lipschitziennes, ce dont on peut s'assurer en se restreignant
à la classe de domaines suivante :
\begin{definition}\label{lbreg}
Un domaine pseudoconvexe borné $\om$ de $\C^k$ est dit $LB$-régulier si $\om$
est $B$-régulier et si pour tout compact $K$ de  $b\om$, il existe une
fonction localement bornée $C_K\geq 0$ et, pour tout $p\in b\om\priv K$ une
fonction $\psi_p$ telle que :
\begin{itemize}
\item $\psi_p\in PSH(\om)\cap\mathcal C^0(\overline{\om})$ et $\psi_p<1$
  sur $\overline{\om}\priv\{p\}$,
\item $\psi_p(p)=1$ et $\psi_p\leq 0$ sur $K$,
\item $\psi_p$ est $C_K(p)$-lipschitzienne en $p$.
\end{itemize}
\end{definition}
\noindent \textbf{Exemples :}
\begin{itemize}
\item[\sbull] Les domaines strictement géométriquement convexes sont $LB$-réguliers.
\item[\sbull] Les domaines de type fini de $\C^2$ sont $LB$-réguliers (\cite{FS}).
\end{itemize}
\begin{rque}\label{pda} 
Le bord d'un domaine $LB$-régulier ne contient pas de disque analytique.
\end{rque}
\begin{prop}\label{antipic}
Soit $\om\subset \C^k$ un domaine $LB$-régulier. 
Pour $U,V$ ouverts de $b\om$ tels que $\overline{U}\cap \overline{V}=\emptyset$, il
existe une constante $C=C(U,V)$ et une fonction $\chi_{U,V}\in PSH(\om)\cap
\mathcal C^0(\overline{\om})$ telle que : 
\begin{description}
\item[i)] $\chi_{U,V}=-1$ sur $U$
\item[ii)] $\chi_{U,V}=0$ sur $V$
\item[iii)] $\chi_{U,V}\leq 0$ sur $\om$
\item[iv)] $\chi_{U,V}$ est $C$-lipschitzienne sur $V$.
\end{description}
\end{prop}
\noindent \underline{Preuve :} Soient $U_1$ tel que $U\Subset U_1\Subset
b\om\priv V$ et $C:=\sup\{C_{\overline{U}_1}(p),\; p\in V\}$ où
$C_{\overline{U}_1}$ est la fonction donnée par la définition \ref{lbreg}. Soit
  $u\in \mathcal C^0(b\om)$  
telle que $-1\leq u\leq 0$, $u=0$ sur $b\om\priv U_1$, $u=-1$ sur $U$. 
Un domaine $LB$-régulier étant en particulier $B$-régulier au sens de
Sibony, il existe $\chi\in PSH(\om)\cap \mathcal
C^0(\overline{\om})$ telle que $\chi_{|b\om}=u$ (auquel cas $\chi$ vérifie
évidemment i), ii) et iii)) (voir \cite{sibony}).
De plus, $\chi$ est donnée par :
$$
\chi=\sup\{v\in PSH(\om)\cap \mathcal C^0(\overline{\om}),\; v_{|b\om}\leq u\}.
$$   
Montrons que $\chi$ vérifie iv). Soit $p\in V$ et $\psi_p$ une fonction \psh\
$C$-Lipschitzienne associée à $p,\overline{U_1}$ par la définition \ref{lbreg}. Par
construction, $u\geq (\psi_p-1)_{|b\om}$, donc $0\geq\chi\geq\psi_p-1$. 
Pour $z\in \om$,
$$
|\chi(z)-\chi(p)|= |\chi(z)|\leq |\psi_p(z)-1|=|\psi_p(z)-\psi_p(p)|\leq C\Vert
z-p\Vert
$$
et $\chi$ est bien $C$-lipschitzienne en $p$.\cqfd
\subsection{Preuve du théorème 2.}
Lorsque le bord contient au moins deux points de faible pseudoconvexité, les
dynamiques sur $\om$ et sur $\spc(b\om)$ sont fortement corrélées :\vspace{0.2cm}\\
\textbf{Fait :} \textit{Sous les hypothèses du théorème 2, considérons une
  sous-suite $f^{n_k}$ tendant vers $a\in b\om$. Alors $f^{n_k}_{|\spc(b\om)}\lra
  a$.} \vspace{0.2cm}

Voyons tout d'abord comment conclure la preuve du théorème 2 à partir de ce Fait.
Comme les limites de
$(f^n)_n$ sont des points de $b\om$ (remarque \ref{pda}), on a l'alternative :
\vspace{0.1cm}

\sbull Soit une sous-suite tend vers $a\in \spc(b\om)$ sur $\om$. Comme
$f^{n_k+1}(z)=f^{n_k}(f(z))$ converge à la fois vers $f(a)$ et $a$ pour $z\in
\om$, le point $a$ est fixé par $f$. Le Fait assure par contre
que l'orbite de $a$ accumule toutes les limites de $(f^n)_n$ puisque $a\in
\spc(b\om)$. Donc 
\begin{equation}\label{thm21}
\text{$(f^n)_n$ converge vers $a$ sur $\om$ et sur $\spc(b\om)$.}
\end{equation}
Pour tout point non-errant $p\in \spc(b\om)$, il existe une suite $n_k\in \N$
  telle que $f^{n_k}(p)$ tend vers $p$. D'après (\ref{thm21}), $p=a$ et donc 
$NW(f_{|b\om})\subset \fpc(b\om)\cup\{a\}$. \vspace{0.1cm}

\sbull Soit toutes les valeurs d'adhérence sont des points de
$\fpc(b\om)$. Soit $a\in b\om$ adhérent à l'orbite d'un point
$p\in \spc(b\om)$.  Quitte à procéder à une double extraction, on peut 
supposer que $f^{n_k}(p)$ tend vers $a$ et que $f^{n_k}$ converge sur
$\om$. Le Fait montre alors que $f^{n_k}$ tend vers $a$ donc que $a\in
\fpc(b\om)$. Ceci montre que  $NW(f_{|b\om})\subset
\fpc(b\om)$.  \vspace{0.2cm}\cqfd
\underline{Preuve du Fait :}
Raisonnons par l'absurde. Soit $p\in\spc(b\om)$ tel que $f^{n_k}(p)$ ne
tend pas vers $a$. Soit $d>0$ tel que $B:=B^{\text{\tiny CR}}(p,d)\Subset
\spc(b\om)$. On peut supposer que $(f^{n_k}(p))_k$ reste dans un ouvert $V$ de
$b\om$ qui n'adhère pas à $a$. Montrons que $b\om\priv\{a\}\subset
\spc(b\om)$. En vertu de la proposition \ref{strucbord},
il suffit d'établir l'existence d'un entier $n$ tel que  
$f^n(B)\supset V$
puisque $V$ est un ouvert arbitrairement gros de $b\om\priv \{a\}$.
Considérons un voisinage $U$ de $a$ tel que
$\overline{U}\cap\overline{V}=\emptyset$ et $\chi_{U,V}$ la fonction donnée
par la proposition \ref{antipic}. 

Comme $(f^{n_k})_k$ converge vers $a$ sur $\om$, il existe $n_\eps\in\N$ tel
que 
$$
f^{n_\eps}\big( \bigcup_{q\in B}\mathcal
C^\frac{\pi}{2}(q_\eps,q)\big)\subset \{\chi_{U,V}<-\frac{1}{2}\}.
$$
D'après la proposition \ref{norm2}, $n_q(f^{n_\eps})\gtrsim \eps^{-1}$ pour
tout $q\in B$ tel que $f^{n_\eps}(q)\in V$. Comme $B\Subset \spc(b\om)$, la
propriété \ref{normcomp} montre donc :
$$
f^{n_\eps}(q)\in V \Longrightarrow \Vert{f_n}'(q)u\Vert\gtrsim
\eps^{-\frac{1}{2}} \hspace{0.5cm}\text{ sur } \{(q,u),\; q\in B,\; u\in
T^\C_qb\om\}. 
$$
Supposons qu'il existe $y\in bf^{n_\eps}(B)\cap V$. Soit  $\gamma:[0,1]\mapsto
V$ un chemin complexe  entre $f^n(p)$ et $y$ de
longueur inférieure à diam$^{\text{\tiny CR}}\,V<+\infty$ (car $\overline{V}$ est
une hypersurface compacte minimale). Comme l'application $f^{n_\eps}$ est un
difféomorphisme $\text{CR}$  local en tout point de $B$, elle permet de
relever $\gamma$ en 
un chemin complexe $\wdt \gamma$ tant que $\wdt \gamma$ ne sort pas de $B$. Le
chemin complexe $\wdt \gamma$ est alors défini sur $[0,\wdt l]$ avec $\wdt l\leq 1$, à
valeurs dans $\overline{B}$ et de longueur $\ell(\wdt \gamma)\geq d$. Comme
$f^{n_\eps}\circ \wdt \gamma=\gamma$ est inclus dans $V$, on a :
$$
\text{diam}^{\text{\tiny CR}} V\geq \ell(\gamma)\gtrsim\ell(\wdt
\gamma)\eps^{-\frac{1}{2}} \approx \eps^{-\frac{1}{2}}.
$$
C'est impossible, donc $bf^{n_\eps}(B)\cap V=\emptyset$. Comme
$f^{n_\eps}(p)\in V$, ceci établit que $f^{n_\eps}(B)\supset V$.\cqfd

\begin{rque}\label{remplacement}
Dans l'énoncé du théorème 1, on peut remplacer l'hypothèse d'existence d'une
fonction \psh\ définissante globale par une hypothèse de $LB$-régularité.
\end{rque}
\section{Auto-applications holomorphes propres des domaines
  disqués.}\label{partiedisque} 
L'objet de cette partie est d'établir le théorème 3. Commençons par
donner les grandes lignes de la démonstration. Compte tenu du théorème
1, on peut supposer que la dynamique de $f$ est non-récurrente.
Nous procédons alors en trois étapes. Dans la première, nous
montrons que le lieu de faible pseudoconvexité est fixé par $f$ et en
déduisons que la dynamique non-récurrente de $f$ est produite
par un point d'attraction $a\in \fpc(b\om)$. Nous montrons dans la deuxième 
étape que le lieu de faible pseudoconvexité est réduit au cercle $C_a$ passant
par $a$. La dernière  
étape exhibe une contradiction à partir de considérations entropiques : 
l'entropie topologique de $f$ est supérieure au logarithme de son degré
et se concentre sur $C_a$, ceci s'avère inpossible lorsque $f$ branche.
\subsection{Fibration de Hopf.}\label{fibrationhopf}
Le bord d'un domaine $\om\Subset \C^2$ disqué et $LB$-régulier a une structure
de fibré en cercles particulièrement utile.

Pour tout $\eta=(z_1,z_2)\in b\om$, nous noterons $D_\eta$ le disque défini par 
$D_\eta:=\{\zeta \eta,\; \zeta\in \D \}$, $C_\eta:=b\D_\eta$ son bord et
$R(\eta)$ son rayon. Rappelons que $\om$ est disqué si et seulement si
$D_\eta\subset \om$ pour tout $\eta\in b\om$. Bien entendu, $C_\eta\subset
b\om$.   
\begin{lemme}
Il existe  un homéomorphisme $h:b\om\lra S^3$ défini par $h(\eta)=\eta/R(\eta)$
commutant aux actions de $\S^1$ :
$$
\forall \zeta\in S^1,\; \forall \eta\in b\om, \; h(\zeta
\eta)=\zeta h(\eta).
$$
\end{lemme}
\noindent \underline{Preuve :} Seule l'affirmation que $h$ est un homéomorphisme est
non triviale. Il suffit de montrer que $h$ est bijective pour l'établir. La
surjectivité est claire car $\om$ est borné. Pour l'injectivité, observons que
deux points distincts $\eta_1,\eta_2$ de $b\om$ tels que $h(\eta_1)=h(\eta_2)$
vérifient $\eta_1=\rho\eta_2$, avec $|\rho|>1$ (ou $|\rho|<1$). Comme $\om$ est disqué,
 si de tels points existaient, la couronne $\{\eta=t\eta_2,\;
 1<|t|<|\rho|\}$ serait incluse dans $b\om$, ce qui contredirait 
la minimalité et donc la $LB$-régularité de $\om$.\cqfd 

\noindent Cet homéomorphisme transporte la fibration de Hopf de $S^3$ sur une
fibration en cercles sur $b\om$ (dont les fibres sont les $C_\eta$, $\eta\in
b\om$) que nous appellerons fibration de Hopf sur $b\om$. On a donc le
diagramme commutatif :
$$
\xymatrix{
    b\om \ar[r]^h \ar[rd]_\pi & S^3 \ar[d]^{\pi_{S^3}} \\
                       & \P^1
  }
$$
On peut ainsi ramener toutes les propriétés de la fibration de Hopf de $S^3$ à
celle de $b\om$. En particulier,  
\begin{prop}\label{noeuds}
Soient $\eta_1,\eta_2\in b\om$ avec $C_{\eta_1}\neq C_{\eta_2}$. Alors
$C_{\eta_1}$ et $C_{\eta_2}$ sont noués dans $b\om$, autrement dit $C_{\eta_1}$ n'est
pas contractile dans $b\om\priv C_{\eta_2}$ (voir \cite{berger},8.6 et 9.4.2).
\end{prop}
Tout chemin noué à un cercle $C_\eta$ et inclus dans un ouvert contractile $U$
de $b\om$ y est noué à $C_a\cap U$. Précisément, 
\begin{lemme}\label{trivloc}
Soit $\pi:b\om\lra \P^1$ la fibration de Hopf et
$\Phi:\D\times]-\eps,\eps[\lra U\subset b\om$ un difféomorphisme fibré
(\textit{i.e.} $\pi\circ\Phi(z,\cdot)$ est constante).
Soient $X\subset U$ et $\gamma$ un lacet dans $U$ dont la projection $\pi\circ
\gamma$ est contractile dans $\D\priv \pi (X)$. Alors $\gamma$ 
est contractile dans $U\priv X$.
\end{lemme}
\noindent \underline{Preuve :}
Ecrivons $\gamma(t)=\Phi(x(t),y(t))$ avec $t\in [0,1]$, $x(t)=\pi\circ\gamma(t)\in \D$ et
$y(t)\in [-\eps,\eps]$.\\
Si  $\pi\circ \gamma$ est contractile, il existe une homotopie $x_s(t)$ entre
$x(t)$ et un point $x$ de $\D\priv \pi(X)$ dans $\D\priv \pi(X)$. Alors
$\gamma_s(t)=\Phi(x_s(t),(1-s)y(t))$ définit une homotopie entre $\gamma(t)$
et $(x,0)$ dans $U\priv X$. \cqfd 

\subsection{Deux lemmes.}
Nous montrons qu'un cercle faiblement pseudoconvexe de la fibration
dont l'image par une application holomorphe propre rencontre le lieu de
stricte pseudoconvexité est isolé dans l'ensemble de faible pseudoconvexité de
$b\om$. 
\begin{lemme}\label{petitbranchement}
Soit $\om\Subset \C^2$ un domaine pseudoconvexe disqué à bord lisse. Soit
$f:\om\lra \om$ une application holomorphe propre, lisse sur
$\overline{\om}$. Si $\eta\in \fpc(b\om)$ et $f(\eta)\in \spc(b\om)$ 
alors  $C_\eta$ est isolé dans l'ensemble des points de faible
pseudoconvexité, en particulier $\fpc(\eta)=C_\eta$.  
\end{lemme}
\noindent \underline{Preuve :}
Raisonnons par l'absurde. Soit $\eta \in \fpc(b\om)$ tel que $f(\eta)\in
\spc(b\om)$. Supposons que $C_\eta$ ne soit pas isolé dans l'ensemble des
points de faible pseudoconvexité. Soit alors une suite de points
$\eta_n$ de $\fpc(b\om)$ tendant vers $\eta$, tels que les cercles $C_{\eta_i}$
soient distincts. Comme $\spc(b\om)$ est ouvert et que $f_{|b\om}$ est
continue, il existe des  voisinages $\Gamma_n$ de $\eta_n$ dans 
$C_{\eta_n}$ tels que $f(\Gamma_n)\subset \spc(b\om)$. Alors le jacobien de
$f$, Jac$(f)$, est nul sur $\Gamma_n$ pour $n$ assez grand (lemme
\ref{strucbord}).  D'après le théorème de Fatou,  
Jac$(f)=0$ sur $\D_{\eta_n}$ pour $n$ grand. L'hypersurface analytique Jac$(f)=0$ 
contient donc une infinité de droites passant par l'origine, c'est impossible.\cqfd 

Le lemme suivant montre que l'image par une certaine itérée de tout cercle de
$b\om$ rencontre un voisinage prescrit d'une valeur d'adhérence de $(f^n)_n$
dans $b\om$. 
\begin{lemme}\label{attradisk}
Soit $\om\Subset \C^2$ un domaine pseudoconvexe disqué et $LB$-régulier. Soit
$f:\om~\lra~\om$ une application holomorphe propre lisse 
sur $\overline{\om}$ dont la
dynamique est non-récurrente. Soit $a\in b\om$ une valeur d'adhérence de
$(f^n)_n$ et $V$ un voisinage de $a$ dans $b\om$. Alors
il existe un entier $n_0$ tel que : 
$$
f^{n_0}(C_\eta)\cap V\neq \emptyset \hspace{0.5cm}\forall \eta\in
b\om.
$$
\end{lemme}
\noindent \underline{Preuve :} 
Comme $\om$ est $B$-régulier, il existe une fonction $u\in C^0(\overline
\om)\cap PSH(\om)$ telle que  $u(a)=1$, $u<1$  sur  $\overline{\om}\priv \{a\}$, 
et  $u=0$  sur $b\om \priv V$. Par hypothèse, il existe un entier $n_0$ tel
que $u\circ f^{n_0}(0)>0$. Le principe du maximum appliqué à $u\circ
f^{n_0}_{|D_\eta}$ montre alors que $f^{n_0}(C_\eta)$ rencontre $V$.\cqfd

\subsection{Démonstration du théorème 3.}
La fonction de jauge est \psh\ et définissante globale. L'application $f$ se
prolonge donc différentiablement sur $\overline{\om}$.
On peut supposer que $\fpc(b\om)$ n'est pas vide et, grâce au théorème 1, que
la dynamique de $f$ est non-récurrente. D'après le théorème 2, il existe alors
une sous-suite $(f^{n_k})_k$ qui converge localement uniformément vers un
point $a\in b\om$  sur $\spc(b\om)$. \vspace{0.2cm}\\ 
\noindent\underline{Etape 1 : $f(\fpc(b\om))= \fpc(b\om)$.} Il suffit de
prouver une inclusion car $f(\spc(b\om))\subset \spc(b\om)$ et $f$ est
surjective. Procédons par l'absurde. Soit $\eta\in \fpc(b\om)$ tel que
$f(\eta)\in \spc(b\om)$. D'après le lemme \ref{petitbranchement},
$\fpc(\eta)=C_\eta$.\\ 
\sbull Montrons qu'il existe $\eta_1$ et $\eta_2$ tels que
$\eta_1$, $\eta_2$, $\eta$ soient deux à deux distincts et 
$$
\left\{
\begin{array}{l}
f(C_{\eta_2})\subset C_{\eta_1}\\
f(C_{\eta_1})\subset C_{\eta}
\end{array}\right.
$$
Pour cela, nous utiliserons le lemme suivant :
\begin{lemme}\label{loc1} Pour tout  voisinage compact $\Gamma$ de $\eta$
  dans $C_\eta$, il existe un point
  $\eta_1$ de $f^{-1}(\Gamma)$ tel que $\Lambda_1(f^{-1}(\Gamma)\cap
  B(\eta_1,\eps))>0$ pour tout $\eps$.   
\end{lemme}
\noindent \underline{Preuve :} Supposons au contraire que pour tout point
$y$ de $f^{-1}(\Gamma)$, il existe un réel $\eps(y)>0$ tel que 
$\Lambda_1\big(f^{-1}(\Gamma)\cap B(y,\eps(y))\big)=0$. Comme $f^{-1}(\Gamma)$
est compact dans $b\om$, on a 
$$
f^{-1}(\Gamma) =\bigcup_{y\in f^{-1}(\Gamma)} B(y,\eps(y))\cap
f^{-1}(\Gamma)=\bigcup_{i=1}^n B(y_i,\eps(y_i))\cap f^{-1}(\Gamma) 
$$
et donc  $\Lambda_1(f^{-1}(\Gamma))\leq \sum_1^n\Lambda_1\big(f^{-1}(\Gamma)\cap
B(y_i,\eps(y_i))\big)=0$. Puis, $f$ étant lisse, $\Lambda_1(\Gamma)=
\Lambda_1\big(f(f^{-1}(\Gamma))\big)=0$ ce qui est absurde. \cqfd 

Soit $\Gamma$ un voisinage compact de $\eta$ dans $C_\eta$ tel que
$f(\Gamma)\subset \spc(b\om)$ et $\eta_1$ donné par le lemme \ref{loc1}. Puisque
$f(\eta_1)\in\fpc(b\om)$ on a aussi $\eta_1\in \fpc(b\om)$. Alors, comme
$f^2(\eta_1)\in f(\Gamma)\subset \spc(b\om)$, le lemme \ref{petitbranchement}
montre que $C_{\eta_1}$ est isolé dans $\fpc(b\om)$. Si $\eps$ est
suffisamment petit, on a donc  
$$
\fpc(b\om)\cap B(\eta_1,\eps)\subset C_{\eta_1}.
$$ 
Comme par ailleurs $f^{-1}(\Gamma)\subset \fpc(b\om)$, on obtient :
$$
\begin{array}{lcl}
0<\Lambda_1(f^{-1}(\Gamma)\cap B(\eta_1,\eps))&=&\Lambda_1\big(f^{-1}(\Gamma)\cap
\fpc(b\om)\cap B(\eta_1,\eps)\big)\\ 
 & \leq & \Lambda_1(f^{-1}(\Gamma)\cap C_{\eta_1})
\end{array}
$$  
Ainsi $\Lambda_1(f^{-1}(\Gamma)\cap C_{\eta_1})>0$, et d'après le théorème de Fatou,
 $f(C_{\eta_1})\subset C_\eta$. On a $C_{\eta_1}\neq C_\eta$ car sinon
 $f(C_\eta)\subset C_\eta\subset \fpc(b\om)$.

De la même façon, on trouve $\eta_2$ ayant les propriétés voulues.\\
\sbull Le théorème de Fatou montre en fait que  
$$
\left\{
\begin{array}{l}
f(\D_{\eta_2})\subset \D_{\eta_1}\\
f(\D_{\eta_1})\subset \D_{\eta}
\end{array}\right.
$$
et il s'ensuit que 
$\{f(0)\}=f(\D_{\eta_2}\cap\D_{\eta_1})\subset\D_{\eta_1}\cap\D_{\eta}=\{0\}$.
Ainsi $f(0)=0$, c'est la contradiction attendue car la dynamique de $f$ est non
récurrente.\vspace{0.2cm}

A ce stade de la preuve, il est bon de remarquer que $a$ ne peut être un point
de stricte pseudoconvexité. En effet, si tel était le cas, les images
successives des cercles de faible pseudoconvexité resteraient en dehors d'un
voisinage fixé de $a$, ce qui est impossible en vertu du lemme
\ref{attradisk}. Le point $a$ est donc nécessairement un point de faible
\vspace{0.2cm} pseudoconvexité. \\
\underline{Etape 2 : $\fpc(b\om)=C_a$.} Montrons tout d'abord que $\fpc(b\om)$
est connexe, il s'agit de voir que $\fpc(\eta)=\fpc(a)$ pour tout $\eta\in
\fpc(b\om)$. Soit $V$ un voisinage de $a$ dans $b\om$, d'après le lemme
\ref{attradisk}, il existe un entier $n_0$ tel que :
$$
f^{n_0}(C_x)\cap V \neq \emptyset \hspace{0.5cm}\forall x\in b\om.
$$
Comme $f^{-1}\big(\fpc(b\om)\big)=\fpc(b\om)$, il existe $\eta_1\in \fpc(b\om)$ tel
que $f^{n_0}(\eta_1)=\eta$. Alors $\fpc(\eta)$ rencontre $V$ car
$\fpc(\eta)\supset f^{n_0}(\fpc(\eta_1))\supset f^{n_0}(C_{\eta_1})$. On a
donc $\fpc(\eta)=\fpc(a)$ car $V$ est un voisinage arbitraire de $a$. 

Montrons maintenant que  $\fpc(b\om)=C_a$.
Notons $X:=\fpc(b\om)=\fpc(a)$. Il s'agit de montrer que $\pi(X)=\pi(a)$. Soit $U$ un
voisinage de $a$ sur lequel il existe une trivialisation locale de la
fibration de Hopf vérifiant : 
$$
\left\{
\begin{array}{l}
\Phi:U \tilde \lra \D\times [-\eps,\eps] \text{ avec }
\pi=\pi_1\circ \Phi\\ 
\Phi(a)=(0,0)\\
\end{array}\right.
$$
Comme $U$ peut être choisi arbitrairement petit, il suffit de montrer que
$\pi(X)$ ne rencontre pas le bord de $\D$. 
Choisissons $C_\eta$ un cercle de $\spc(b\om)$ et $n_0$ un entier tel que
$f^{n_0}(C_\eta)\subset U$ (voir le théorème 2).  
Le cercle $C_\eta$ n'est pas contractile dans $b\om \priv X$ car $X$ contient des
cercles de la fibration de Hopf (propriété \ref{noeuds}).  
Comme, d'après la première étape,  $f(\fpc(b\om))\subset \fpc(b\om)$, $f^{n_0}$
induit un revêtement fini de $\spc(b\om)$ sur lui-même (voir
proposition \ref{revetbord}). Il s'ensuit que $f^{n_0}(C_\eta)$ n'est pas contractile dans
$U\priv X$, et donc que  $\pi(f^n(C_\eta))$ ne l'est pas dans $\D\priv \pi(X)$
(lemme \ref{trivloc}). Comme $\D$ est assimilable à un ouvert de $\R^2$, cela
n'est possible que si $\pi(X)\cap b\D=\emptyset$.\vspace{0.2cm}\\
\underline{Etape 3 : $f$ est injective.} Soit $d$ le degré
topologique de $f_{|b\om}$; c'est aussi celui de  $f$. D'après  un théorème de
Misiurevitcz-Przytycki (\cite{misprz}), l'entropie topologique de $f_{|b\om}$
est minorée :  
$$
h_{top}(f_{|b\om})\geq \log d.
$$
D'après le principe variationnel, cette entropie est portée par l'ensemble
non-errant (\cite{katok}, formule 3.3.1) :
$$
h_{top}(f_{|b\om})=h_{top}(f_{|NW(f_{|b\om})}).
$$
Or d'après le théorème 2 et les deux premières étapes, on a
$NW(f_{|b\om})\subset C_a$ si bien que 
$$
h_{top}(f_{|C_a})\geq \log d.
$$
Comme $f(C_a)=C_a$, on voit facilement que $f(\D_a)=\D_a$. Alors $f_{|\D_a}$
est un produit de Blaschke fini qui, vu comme fraction  rationnelle $\wdt f$,
est de degré $\wdt d=deg(f_{|\D_a})$. Or, comme l'a montré Gromov (\cite{gromov}), 
$h_{top}(\wdt f)=\log \wdt d$, et donc
$$
h_{top}(f_{|C_a})=h_{top}(\wdt f_{|C_a})\leq \log\wdt d.
$$  
Alors $deg(f_{|\D_a})=\wdt d\geq d=deg f$ ce qui montre que l'ensemble des
valeurs critiques $f(V_f)$ ne
contient pas $\D_a$. Or, il résulte des deux premières étapes que $V_f\subset
\D_a$ et $f(\D_a)=\D_a$. Ainsi, $V_f=\emptyset$ et $f$ est injective.\cqfd


\appendix
\section{Existence de limite normale de $\rho$ en un point de stricte
  pseudoconvexité.}\label{limnormro}
Le but de cet appendice est de montrer le lemme \ref{limnorm}. Nous allons 
 prouver que l'ensemble des points en lesquels $\rho$ possède des limites
radiales dans $M$ est dense dans $\spc(b\om)$. Fixons pour cela un ouvert
arbitraire $U\subset \spc(b\om)$.
 
\begin{figure}[h]
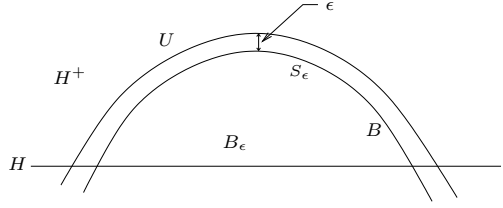

\begin{center}
\input dessin2.pstex_t
\end{center}
\caption{Notations de la preuve du lemme \ref{limnorm}.}
\label{sitrud}
\end{figure}

Commençons par introduire les notations utilisées dans la preuve (on se
reportera avec avantage à la figure \ref{sitrud} qui les illustre).\\
Tout d'abord, quitte à restreindre $U$ et à effectuer un changement de variable
holomorphe, on peut supposer que $U$ est formé de points de stricte convexité
de $b\om$.
Soit alors $H$ un hyperplan réel de $\C^k$ tel que $H\cap U\Subset U$. Soit $H^+$
le demi-espace délimité par $H$ tel que $H^+\cap U\Subset U$. Pour  $\eps>0$,
$S_\eps$ désignera l'ensemble des points de $\om$ à distance $\eps$ de $U$ :   
$$
\begin{array}{l}
S_\eps=\{p+\eps\vec N(p),\; p\in U\}\cap H^+\\
B_\eps=\{p+r\vec N(p),\; p\in U,r>\eps\}\cap H^+=\cup_{r>\eps} S_\eps ,
\end{array}
$$ 
enfin, $B=\displaystyle \bigcup_{\eps>0} B_\eps$.

Nous utiliserons les lemmes suivants :
\begin{lemme}\label{maxan}
Soit $A$ un ensemble analytique de codimension 1 de $\C^k$ et $\Lambda$ une
fonction \psh\  continue de $\C^k$, non constante sur $A$. 
Alors $\Lambda_{|A}$ n'a pas de maximum local.  
\end{lemme}
\begin{lemme}\label{radpropre}
Soit $\Phi:\D\lra \D$ holomorphe, radialement propre. Alors
Im$\Phi$ est dense dans $\D$. 
\end{lemme}
\begin{lemme}\label{hyp}
Soit $M$ une surface de Riemann de genre fini. Il existe un compact $K$ de
$M$ tel que toute composante connexe de $M\priv K$ est
\begin{itemize}
\item soit un anneau $\mathcal A(r,1)$ avec $\S_r\subset K$,
\item soit un disque épointé $\D^*$ avec $\S_1\subset K$. 
\end{itemize}
\end{lemme}
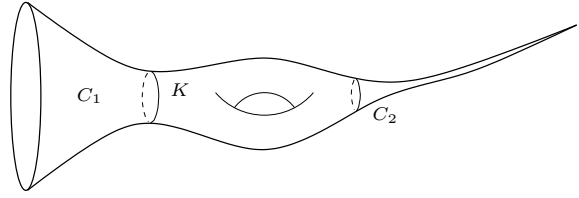
\begin{figure}[h]
\begin{center}
\input hyperbolic3.pstex_t
\end{center}
\end{figure}

\paragraph{} Nous procédons par l'absurde. On montre par un argument de
catégories de Baire que, quitte à restreindre $B$, l'image de $B$ par  
$\rho$ évite un compact $K$ de $M$ préalablement fixé.
En restreignant $\rho$ à une droite complexe transverse à $b\om$,
on obtient une application radialement propre d'un disque dans une composante
de $M\priv  K$. Cette situation se révèle être impossible.\vspace{0.2cm}\\ 
\underline{Preuve du lemme \ref{limnorm}:} Supposons que  $\rho^*(p)\cap
M=\emptyset$ pour tout $p\in U$, \textit{i.e.} $\rho$ est radialement propre
en les points de  $U$. Soit $K$ un compact de $M$ donné par le lemme 
\ref{hyp}. 

Montrons dans un premier temps que, quitte à restreindre $U$ :
\begin{equation}\label{a1}
\rho(B)\cap K=\emptyset
\end{equation}
Pour cela, soit $r_j\searrow 0$, et  
$$
E_{i}=\left\{ p\in U|\; \forall j\geq i,\;\rho(p+r_j\vec
  N(p))\notin K\right\}.
$$ 
Comme $\rho^*(p)\cap M=\emptyset$ pour tout $p\in U$ on a  $\bigcup_{i}
E_{i}=U$. 
D'après le théorème de Baire, il existe donc $j_0\in \N$ tels que $E_{j_0}$ soit
d'intérieur non vide. Quitte à restreindre $U$, on peut supposer que
$E_{j_0}=U$. Alors 
\begin{equation}\label{a2}
\rho(S_{r_j})\cap K=\emptyset \hspace{0.5cm}\forall j\geq j_0.
\end{equation}
Soit $z\in B$ et $\F$ la composante connexe de $\rho^{-1}(\rho(z))$ qui
contient $z$. C'est un
ensemble analytique de codimension 1. Soit $\Lambda$ une forme linéaire telle 
que $H^+=\{\Lambda>0\}$. Le principe du maximum appliqué à $\rho_{|\F}$ montre
que $\F\cap S_{r_j}\neq \emptyset$, donc que $\rho(B_{r_j})\subset \rho(S_{r_j})$.
Grâce à (\ref{a2}), ceci termine la preuve de (\ref{a1}) puisque
$$
\rho(B)\cap K=\rho(\bigcup_{j\geq j_0} B_{r_j})\cap
K =\bigcup_{j\geq j_0}(\rho(B_{r_j})\cap K)=\emptyset.
$$

Dans un second temps, nous restreignons $\rho$ à une droite complexe
générique. Ceci produit encore une application propre par un argument de
Lindelöff-Cirka. \vspace{0.2cm}\\  
\noindent\textbf{Fait}\textit{
Soit $l$ une droite complexe sur laquelle $\rho$ est non constante,
qui intersecte $B$ et telle que $l\cap bB\subset U$. Alors $\rho_{|l\cap B}$
est radialement propre.}\vspace{0.2cm}

Concluons la preuve à partir de ce Fait.
Comme $B$ est connexe, $\rho(B)$ est inclus dans
une composante connexe de $M\priv K$ d'après (\ref{a1}). L'application 
$\rho'$ est donc à valeurs dans  un disque épointé ou un 
anneau $\mathcal A(r,1)$. Le Fait précédent montre que les limites radiales de
$\rho'$ sont nulles dans le premier cas - c'est impossible d'après le théorème
de Fatou - ou de module $1$ dans le second.  
On peut donc voir $\rho'$ comme une fonction holomorphe à valeurs dans $\D$
radialement propre et évitant un ouvert de $\D$. En composant à gauche 
$\rho':B'\lra \D$ par une représentation conforme de $\D$ sur $B'$, on produit
une fonction holomorphe de $\D$ dans $\D$ radialement propre dont l'image
n'est pas dense. Ceci est
impossible d'après le lemme \ref{radpropre}.\cqfd
\underline{Preuve du Fait :}
Par commodité, nous noterons $B'=l\cap B$, $U'=l\cap U= bB'$,
et $\rho'=\rho_{|B'}$. Pour $p\in U'$, définissons aussi $\vec{N}_l(p)$ le
vecteur unitaire normal  à $U'$ rentrant dans $B'$. 
Montrons que la propreté radiale de $\rho$ implique celle de $\rho'$.  

Comme $i\vec{N}_l(p)\in T_p(U\cap l)$,
il s'agit d'un vecteur  de la forme $ai\vec{N}(p)+y$ avec $y\in
T^\C_pU$. On en déduit que 
$$
\vec{N}_l(p)=a\vec{N}(p)+x,\hspace{0.5cm} x\in T_p^\C U
$$
Remarquons que $a$ est positif car $\vec{N}(p)$ et $\vec N_l(p)$ sont tous
deux des vecteurs rentrants. 

Supposons alors qu'il existe une suite de réels $t_k \searrow 0$
telle que $\rho'(p+t_k\vec{N}_l(p))$ converge vers $\beta \in M$. 
Ceci signifie $\rho(p+at_k\vec{N}(p)+t_kx)\lra \beta$.
Des estimations standard de distance de Kobayashi  montrent
que $d_{K_B}(p+at_k\vec{N},p+at_k\vec{N}+t_kx)\in O(\sqrt t_k)$. L'application
$\rho$ étant contractante pour la  
métrique de Kobayashi, $d_{K_M}(\rho(p+at_k\vec{N}),
\rho(p+at_k\vec{N}+t_kx))$ tend vers $0$. Donc $\rho(p+at_k\vec{N})$ tend vers
$\beta\in M$ ce qui est faux puisque $\rho$ est radialement propre.  \cqfd
\bibliography{bib1.bib}
\bibliographystyle{abbrv}

\bigskip
Emmanuel Opshtein\\
Laboratoire Emile Picard \\
UMR 5580, UFR MIG \\
Université Paul Sabatier \\
118, route de Narbonne \\
31062 Toulouse cédex 4 \\
opshtein@picard.ups-tlse.fr

\end{document}